\documentclass[a4paper, 11pt]{amsart}
\pdfoutput=1

\usepackage{amssymb}
\usepackage{amsfonts}
\usepackage{amsthm}
\usepackage{stmaryrd}
\usepackage{kpfonts} 
\usepackage[mathscr]{eucal}

\usepackage[english]{babel}

\usepackage[T1]{fontenc}
\usepackage[utf8]{inputenc}

\usepackage[margin=1in]{geometry}

\usepackage{tikz-cd}
\usepackage{tikz}
\usetikzlibrary[trees]
\usepackage{enumitem}
\usepackage{hyperref}
\usepackage{aliascnt}
\usepackage{adjustbox}
\usepackage{varioref}
\usepackage{xstring}
\usepackage{xcolor}
\usepackage{bbm}

\newtheorem{thmintro}{Theorem}

\newcommand{\inewtheorem}[2]{
	\newaliascnt{#1}{thmintro}
	\newtheorem{#1}[#1]{#2}
	\aliascntresetthe{#1}
}
\inewtheorem{propintro}{Proposition}
\inewtheorem{corintro}{Corollary}

\newtheorem{theorem}{Theorem}
\newcommand{\jnewtheorem}[2]{
	\newaliascnt{#1}{theorem}
	\newtheorem{#1}[#1]{#2}
	\aliascntresetthe{#1}
}
\numberwithin{theorem}{section}
\jnewtheorem{lemma}{Lemma}
\jnewtheorem{proposition}{Proposition}
\jnewtheorem{corollary}{Corollary}
\jnewtheorem{conjecture}{Conjecture}

\theoremstyle{definition}
\jnewtheorem{definition}{Definition}
\jnewtheorem{example}{Example}
\jnewtheorem{remark}{Remark}
\jnewtheorem{question}{Question}
\jnewtheorem{fact}{Fact}
\jnewtheorem{construction}{Construction}
\jnewtheorem{notation}{Notation}
\jnewtheorem{convention}{Convention}

\labelformat{theorem}{Theorem #1}
\labelformat{lemma}{Lemma #1}
\labelformat{proposition}{Proposition #1}
\labelformat{corollary}{Corollary #1}
\labelformat{definition}{Definition #1}
\labelformat{example}{Example #1}
\labelformat{section}{Section #1}
\labelformat{subsection}{Subsection #1}
\labelformat{equation}{(#1)}
\labelformat{remark}{Remark #1}
\labelformat{exercise}{Exercise #1}
\labelformat{question}{Question #1}
\labelformat{chapter}{Chapter #1}
\labelformat{construction}{Construction #1}
\labelformat{conjecture}{Conjecture #1}
\labelformat{notation}{Notation #1}
\labelformat{convention}{Convention #1}

\labelformat{thmintro}{Theorem #1}
\labelformat{corintro}{Corollary #1}
\labelformat{propintro}{Corollary #1}

\numberwithin{table}{subsection}
\labelformat{table}{Table #1}

\hypersetup{
	pdftitle={},
	pdfauthor={Gregoire Marc},
	bookmarksnumbered=true,     
	bookmarksopen=true,         
	bookmarksopenlevel=1,       
	colorlinks=true,   
	linkcolor = {cyan!60!black},   
	urlcolor = {cyan!60!black},   
	citecolor = {cyan!60!black},   
	pdfencoding=unicode,      
	pdfstartview=Fit,           
	pdfpagemode=UseOutlines,    
	pdfpagelayout=OneColumn,
	breaklinks=true,
	linktocpage
}

\def\FF{{\mathbf{F}}}
\def\fG{{\mathbb{G}}}

\def\sF{{\mathscr{F}}}

\def\sO{{\mathscr{O}}}

\def\sS{{\mathscr{S}}}

\def\fG{{\mathfrak{G}}}

\def\cC{{\mathcal{C}}}
\def\cD{{\mathcal{D}}}
\def\cE{{\mathcal{E}}}

\def\cG{{\mathcal{G}}}
\def\cI{{\mathcal{I}}}
\def\cJ{{\mathcal{J}}}

\def\cL{{\mathcal{L}}}
\def\cO{{\mathcal{O}}}

\mathchardef\bigppchar"1403
\def\dobigpp#1{\vcenter{#1\kern.2ex\hbox{$\bigppchar$}\kern.2ex}}
\mathchardef\bigppdomchar"1400
\def\dobigppdom#1{\vcenter{#1\kern.2ex\hbox{$\bigppdomchar$}\kern.2ex}}

\mathchardef\pp"2403 
\mathchardef\ppdom"2400 

\newcommand\PP[1]{\mathscr{P}_{#1}}

\newcommand\Cat{\operatorname{Cat}}

\newcommand\Id[1]{\operatorname{id}_{#1}}

\newcommand\opo[1]{{#1}^{\operatorname{op}}}

\newcommand\fun{\operatorname{Fun}}

\newcommand\Topo{\operatorname{Top}}

\newcommand\Psh{\operatorname{Psh}}
\newcommand\Arr{\operatorname{Arr}}

\newcommand\num[1]{\langle#1\rangle}

\newcommand\suborb[1]{\operatorname{SubOrb}_{#1}}
\newcommand\OrbG{\suborb{\sO_G}}

\newcommand\Ob{\operatorname{Ob}}
\newcommand\func{\operatorname{Fun}}
\newcommand\Hom{\operatorname{Hom}}

\newcommand\OG{\mathscr{O}_G}

\newcommand\catG{\operatorname{Cat}^G}
\newcommand\CatG{\mathbf{Cat}^G}

\newcommand{\gset}{\operatorname{Fin}^G}
\newcommand{\gsett}{\operatorname{Set}^G}
\newcommand\set{\operatorname{Set}}

\newcommand\gope{\operatorname{Op}\left(\gsset\right)}

\newcommand\gOpe{\mathbf{Op}^G}
\newcommand{\Op}{\operatorname{Op}}

\newcommand{\hset}[1]{\operatorname{Fin}^{#1}}

\newcommand\gtop{\operatorname{Top}^G}

\newcommand\Alg{\operatorname{Alg}}
\newcommand\gsset{\operatorname{SSet}^G}
\newcommand\sset{\operatorname{SSet}}
\newcommand\gs{\mathscr{S}^G}
\newcommand\IndG{\operatorname{I}_G}
\newcommand\End{\operatorname{End}}
\newcommand\FFF[1]{\mathbf{F}{#1}}
\newcommand\conss{\mathcal{G}}

\newcommand\SpaG{\mathscr{A}_{\Ind}}

\newcommand\ind[2]{\operatorname{Ind^{#1}_{#2}}}
\newcommand\res[2]{\operatorname{Res}^{#1}_{#2}}
\newcommand\spa[2]{\mathscr{A}(#1,#2)}
\newcommand\Mack{\operatorname{Mack_{\Ind}}}

\newcommand\Cmon{\operatorname{CMon}^G}

\newcommand\NSSMI{\mathcal{I}\mathbf{SMSt}}

\newcommand\NSStMI{\mathcal{I}\mathbf{SMStg}}

\newcommand\NSSlaxI{\mathcal{I}\mathbf{SMLax}}

\newcommand\OAlgL{\mathbf{AlgLax}_{\ma{O}}}
\newcommand\OAlg{\mathbf{Alg}_{\ma{O}}}
\newcommand\OAlgSt{\mathbf{AlgSt}_{\ma{O}}}

\newcommand\OAlgStg{\mathbf{AlgStg}_{\ma{O}}}
\newcommand\EII{\mathcal{E}_{\mathcal{I}}}
\newcommand\Ninf{\mathbf{Alg}_{\EII}}
\newcommand\Bur{\operatorname{A}_G} 
\newcommand\Buri{\operatorname{A}_{\Ind}} 
\newcommand\Sing{\operatorname{Sing}}
\newcommand\SN{F(S_{\Ind})}

\newcommand\nope{N_{\infty}\text{-}\operatorname{Op}}
\newcommand\Nope{N_{\infty}\text{-}\mathbf{Op}}
\newcommand\ma[1]{\mathcal{#1}}

\newcommand\FVI{F(V_{\Ind})}

\newcommand\Mor{\operatorname{Mor}}

\newcommand\TTT{\operatorname{Tr}}

\newcommand\Oalg{\operatorname{Alg}_{\ma{O}}}

\newcommand\EI{\mathcal{SM}_{\Ind}}

\newcommand\Ind{\mathcal{I}}

\title{A higher Mackey functor description of algebras over an $N_{\infty}$-operad}
\author{Gr\'{e}goire Marc}

\begin{document}
	
\maketitle

\begin{quote}
\begin{center}
\textbf{Abstract}
\end{center}
Suppose $G$ is a finite group. In this paper, we construct an equivalence between the $\infty$-category of algebras over an $N_{\infty}$-operad $\ma{O}$ associated with a $G$-indexing system $\Ind$ and the corresponding $\infty$-category of higher incomplete $\Ind$-Mackey functors with value in spaces. We use the universal property of the incomplete $(2,1)$-category of spans of finite $G$-sets $\SpaG$ to construct a functor from $\SpaG$ to the $2$-category of $\Ind$-normed symmetric monoidal categories of Rubin. We then show that the left Kan extension of the composition of this functor with the core functor along the Yoneda embedding is an equivalence.

\end{quote}

\tableofcontents

\section{Introduction}

Let $G$ be a finite group and consider a set $X$ endowed with an action of $G$. For every subgroup $H$ of $G$, we can consider the set of $H$-fixed elements of $X$
$$
X^H=\{x\in X~|~ h\cdot x=x~\text{for all}~ h\in H \}.
$$ 
For every pair of subgroups $K\leqslant H$ of $G$ and every element $g$ of $G$, there is a restriction map $\operatorname{Res}^{H}_{K}:X^H\hookrightarrow  X^K$ (the inclusion of $H$-fixed elements in $K$-fixed elements) and a conjugation map $\operatorname{c}_g:X^H \to X^{gHg^{-1}}$ that sends $x$ to $g\cdot x$. These data assemble into a presheaf $F_X:\opo{\OG}\to \set$ over the category of transitive $G$-equivariant sets $\OG$ that sends the coset $G/H$ to the set $X^H$. This construction forms a fully faithful functor $F_{(-)}:\gsett \to \fun (\opo{\OG},\set)$ that is, in general, not an equivalence of categories. However, if $X$ is a topological space, then, up to weak homotopy equivalence, every topological presheaf $P:\opo{\OG} \to \Topo$ has the form $F_X$ for some $G$-equivariant topological space $X$. Moreover, for the right homotopy theory on $G$-equivariant spaces, \emph{Elmendorf's Theorem} \cite{Elm} states that this functor is an equivalence of homotopy theories. Elmendorf's Theorem can be considered to be one of the cornerstones of equivariant homotopy theory as it brings upon it categorical features, either with the theory of \emph{parametrized $\infty$-categories} developed by Barwick,
Dotto, Glasman, Nardin, and Shah \cite{para, Sha1}, or with the more general theory of \emph{categories internal in an $\infty$-topos} of Martini and Wolf \cite{Mar1,MW2021, Mar2}.

 Suppose now that $X$ is a commutative monoid with a $G$-action compatible with its monoid structure. The presheaf provided by the previous construction comes with an extra structure given by the monoid law. For every pair of subgroups $K\leqslant H$ of $G$, there is a transfer map $\TTT_K^H:X^K \to X^H$ defined by
$$
\TTT_K^H (x)=h_1 \cdot x + \cdots + h_{[H:K]}\cdot x
$$
where $h_1,\ldots, h_{[H:K]}$ are representatives of cosets in $H/K$. These transfers, along with the previous structure of restriction and conjugation maps, assemble into a \emph{Mackey functor} (with value in the category of sets) as defined by Dress in \cite{dress1971notes}. A Mackey functor $\ma{M}_X:\opo{\Bur} \to \set$ is a presheaf over the \emph{Burnside category} $\Bur$ that preserves finite products. The Burnside category $\Bur$ is the category whose objects are finite $G$-equivariant sets and whose morphisms are isomorphism classes of spans between these $G$-equivariant finite sets. Mackey functors play an important role in representation theory, equivariant homotopy theory, and they have many applications, for example in group cohomology and in the decomposition of classifying spaces with the work of Webb \cite{MR1233331}. The construction of a Mackey functor from a commutative monoid with a compatible $G$-action provides a fully faithful functor $\ma{M}_{(-)}:\Cmon \to \func_{\times}(\opo{\Bur},\set)$ that is again not an equivalence of categories in general. One might now be interested in promoting this construction to an equivalence of homotopy theories, in analogy with Elmendorf's Theorem mentioned above. 

 However, unlike classical algebra, in homotopy theory, commutativity is not a property, but a structure. Given a topological space $X$ with an associative and unitary operation $\mu:X\times X\to X$, $X$ is considered as being homotopy commutative if there exists a homotopy $H$ as depicted in the following diagram
$$
\begin{tikzcd}[row sep=small]
X\times X \arrow[dd,"\operatorname{twist}"'] \arrow[rd,"\mu"name=O]& \\
& X \\
|[alias=Z]| X\times X \arrow[ru,"\mu"'] & 
\arrow[from=Z,to=O,Rightarrow,shorten >=3mm,shorten <=4mm,"H"]
\end{tikzcd}
$$
and if higher coherences are satisfied (requiring the data of higher homotopies).
The data of these homotopies are part of the commutative structure of $X$ and can be encoded with May's theory of \emph{$E_{\infty}$-operads} \cite{May}. An $E_{\infty}$-operad is a topological operad $\ma{O}$ such that $\ma{O}_n$ is contractible and has a free $\Sigma_n$-action. This theory of ``commutativity up to homotopy'' is crucial in topology, as most spaces with a multiplication, such as infinite loop spaces, are not strictly commutative.  

  If $X$ is a topological space with a $G$-action and a compatible structure of $E_{\infty}$-algebra (meaning that $X$ is an algebra over an $E_{\infty}$-operad $\ma{O}$ and that all the corresponding operations and homotopies between them are $G$-equivariant), then one can still consider the topological presheaf $F_X:\opo{\OG}\to \Topo$ associated with $X$ but the transfers do not exist anymore, as the operation $\mu$ of $X$ is only commutative up to homotopy. To solve this issue, $E_{\infty}$-operads can be upgraded into an equivariant analogue, the \emph{$N_{\infty}$-operads} introduced by Blumberg and Hill in \cite{BH1}. In addition to the classical operations of an $E_{\infty}$-algebra, an $N_{\infty}$-operad encodes norms that can be thought of as equivariantly twisted versions of these operations. 
Unlike the non-equivariant setting, not all $N_{\infty}$-operads are equivalent, and the homotopy category of $N_{\infty}$-operads is equivalent to the poset of the so-called \emph{$G$-indexing systems}, by the work of Blumberg and Hill \cite{BH1}, Rubin \cite{Rub1}, and Bonventre and Pereira \cite{BonPer1}. A $G$-indexing system $\Ind$ is a collection, for every subgroup $H$ of $G$, of finite $H$-equivariant sets satisfying relations of compatibility, as described in \cite{BH1}. If $\Ind$ is a $G$-indexing system, then every finite $H$-equivariant set $T$ that belongs to $\Ind$ encodes an operation that is, roughly speaking, ``$T$-equivariant''.

Suppose now that $X$ is an algebra over an $N_{\infty}$-operad $\ma{O}$ associated with a $G$-indexing system $\Ind$. If $K \leqslant H$ are subgroups of $G$ such that $H/K$ belongs to $\Ind$, then we may construct a transfer map $\TTT_{K}^{H}:X^K \to X^H$. However, the construction of $\TTT_{K}^{H}$ requires the choice of an $H/K$-equivariant operation on $X$, that is only unique up to homotopy.
The homotopy groups of $X$ thus give rise to \emph{$\Ind$-incomplete Mackey functors}, defined in terms of an incomplete version of the Burnside category $\Buri$, as studied by Blumberg and Hill in \cite{MR3773736}. This structure cannot be na\"ively promoted to the space $X$ itself, as the transfers depend on the choice of homotopies. To take these homotopies into account, the \emph{$\Ind$-Burnside category} $\Buri$ can be upgraded into a $(2,1)$-category, the \emph{$\Ind$-effective Burnside category $\SpaG$}, remembering isomorphisms between spans. The $(2,1)$-category $\SpaG$ leads also to a notion of \emph{higher $\Ind$-incomplete Mackey functors} (with value in the $\infty$-category of spaces) defined as presheaves (of spaces) over $\SpaG$ that preserve finite products. In this paper, we construct an equivalence between the $\infty$-category of (higher) $\Ind$-incomplete Mackey functors and the $\infty$-category of algebras in $G$-topological spaces over an $N_{\infty}$-operad $\ma{O}$ associated with $\Ind$. More precisely, the main result of this paper is the following theorem.  

\begin{thmintro}\label{THETH}(\ref{THETHH})
Let $G$ be a finite group, let $\Ind$ be a $G$-indexing system and denote by $\EII$ Rubin's $N_{\infty}$-operad associated with $\Ind$. The $\infty$-category obtained from $\Alg_{\EII}\left(\gtop\right)$ by inverting $G$-weak homotopy equivalences (i.e., continuous $G$-equivariant maps $f:X\to Y$ such that, for every subgroup $H$ of $G$, $f^H:X^H\to Y^H$ is a weak homotopy equivalence) is equivalent to the $\infty$-category of $\Ind$-incomplete Mackey functors $\Mack=\func_{\times}(\opo{\SpaG},\sS)$ where $\sS$ is the $\infty$-category of spaces and $\SpaG$ is the incomplete $(2,1)$-category of spans of finite $G$-equivariant sets associated with $\Ind$.
\end{thmintro}

\subsection{General motivations and related work}
Our main motivation for \ref{THETH} is to compare the theory of $N_{\infty}$-algebras with the theory of parametrized $\mathscr{T}$-semi-additivity as developed by Cnossen, Lenz, and Linskens in \cite{CLL1}, which generalises the work of Nardin in \cite{dede} and \cite{dede2}. As we already mentioned, one important feature of Elmendorf's Theorem is that the theory of parametrized $\infty$-categories gives a categorical framework for equivariant homotopy theory through the notion of $G$-$\infty$-categories (here we mean $\OG$-parametrized $\infty$-categories). This has been motivated by Guillou and May's \cite[Theorem 1.14]{GM}, which describes genuine equivariant spectra in terms of spectral Mackey functors. The main theorem of the present paper can be thought of as an unstable version of Guillou and May's Theorem. The poset of $G$-indexing systems or, equivalently, the poset of $N_{\infty}$-operads is equivalent to the poset of orbital subcategories of $\OG$ in the sense of \cite{CLL1} and it follows that every $G$-indexing system $\Ind$ corresponds to an associated orbital subcategory $\mathscr{T}$ of $\OG$. Every orbital subcategory $\mathscr{T}$ of $\OG$ comes with an associated notion of $\mathscr{T}$-semi-additivity and it has been shown, first by Nardin in \cite[Theorem 6.5]{dede} and \cite[Theorem 2.32]{dede2} in the genuine case (that is for the maximal $G$-indexing system), and more recently by Cnossen, Lenz, and Linskens in the general case \cite[Corollary 9.14]{CLL3}, that the $G$-$\infty$-category of $\Ind$-incomplete Mackey functors is the universal $\mathscr{T}$-semi-additive $G$-$\infty$-category over equivariant spaces. One of the consequences of \ref{THETH} should be that the $G$-$\infty$-category of $\cO$-algebras for an $N_{\infty}$-operad $\ma{O}$ associated with a $G$-indexing system $\Ind$ is the universal $\mathscr{T}$-semi-additive $G$-$\infty$-category over equivariant spaces where $\mathscr{T}$ is the orbital subcategory of $\OG$ that corresponds to $\Ind$.

The strategy of proof that we adopt to show \ref{THETH} is similar to the techniques used in \cite{TOB1}. Note that it should be possible to establish a formal comparison between $\cO$-algebras for an $N_{\infty}$-operad $\cO$ and incomplete Mackey functors using the theory of parametrized operads \cite{DJ} by showing, for every simplicial $G$-operad $\ma{O}$ such that the $\Sigma_n$-action on $\cO_n$ is free, that the nerve functor $N \left(\Alg_{\ma{O}}\left(\gsset\right)\right)\to \Alg_{N^{\otimes}(\ma{O})}\left(\gs\right)$ of \cite[Theorem IV]{Bon} is an equivalence of $\infty$-categories. However, this would require a comparison, for every operad $\ma{O}$ satisfying these conditions, between free $\ma{O}$-algebras over $G$-equivariant spaces and the corresponding free algebras over the nerve of $\ma{O}$.

\subsection{Structure of the paper}
In \ref{PreNinf}, we first recall Blumberg and Hill's notions of $N_{\infty}$-operads and $G$-indexing systems and the relation between them. We then define the $\infty$-category $\mathbf{Alg}_{\cO}$ of algebras over an $N_{\infty}$-operad $\ma{O}$ and exhibit an explicit description of it. In \ref{Norcat}, we recall Rubin's notion of $\Ind$-normed symmetric monoidal categories for a given $G$-indexing system $\Ind$ and we give an explicit construction of the operad $\EI$ that encodes them. We then show that the nerve of $\EI$, denoted by $\EII$, is an $N_{\infty}$-operad associated with $\Ind$. In \ref{guimay}, we begin by defining the $\Ind$-effective Burnside category $\SpaG$ and the associated $\infty$-category of $\Ind$-Mackey functors $\Mack$. We then construct explicitly the free $\EII$-algebra over a finite $G$-set $A$ as the nerve of an $\Ind$-normed symmetric monoidal groupoid $\sF(A)$ and we use a $2$-categorical unfurling to construct a functor $\theta:\SpaG \to \mathbf{Alg}_{\EII}$ that sends $A$ to the free $\EII$-algebra over $A$. Finally, we use the results of \cite{greg3} to prove \ref{THETH} (\ref{THETHH}) by showing that the left Kan extension of $\theta:\SpaG \to \mathbf{Alg}_{\EII}$ along the Yoneda embedding $\SpaG \to \Mack$ is an equivalence of $\infty$-categories. 
In \ref{ApA}, we show a way to see the $(2,1)$-category of $\Ind$-normed symmetric monoidal groupoids as a full subcategory of $\mathbf{Alg}_{\EII}$. 

\subsection{Notation and conventions} Throughout this paper, $G$ is a finite group and for every category $\cC$, we denote by $\cC^G$ the category of $G$-objects in $\cC$. A \emph{$G$-category} is an internal category to $\gsett$ or, equivalently, a $G$-object in the category of categories. If $\cC$ and $\cD$ are $G$-categories, then we denote by $\func_G (\cC,\cD)$ the category of $G$-functors and $G$-natural transformations between $\cC$ and $\cD$ and by $\func (\cC,\cD)$ the $G$-category of nonequivariant functors and nonequivariant natural transformations. We denote by $\CatG$ the $2$-category of small $G$-categories.
For every closed symmetric monoidal category $\cC$, we denote by $\Op (\cC)$ the category of operads in $\cC$. If $\ma{O}$ is an operad in $\cC$, then we denote by $\ma{O}_n$ its n-th level, that is the associated object of $\cC^{\Sigma_n}$. We denote by $\sS$ the $\infty$-category of spaces.

\subsection{Acknowledgements}
The work presented in this paper is a part of my PhD project. I would like to thank
my PhD supervisor Magdalena Kędziorek for her guidance and various enlightening conversations. I would also
like to thank Miguel Barrero and Niall Taggart for their very helpful comments on an earlier draft of this paper.
During the writing of this paper, the author was supported by an NWO grant Vidi.203.004.
\section{Preliminaries on \texorpdfstring{$N_{\infty}$-operads}{N-infinity operads}}\label{PreNinf}

\subsection{\texorpdfstring{$N_{\infty}$-operads}{N-infinity operads} and \texorpdfstring{$G$}{G}-indexing systems}
In this subsection, we recall Blumberg and Hill's notions of $N_{\infty}$-operads in the setting of $G$-simplicial sets, of $G$-indexing systems, and the relation between them. We refer the reader to \cite{BH1} or \cite{Rub1} for more details on the subject. 

\begin{definition}
	Let $n\ge 0$. A \emph{graph subgroup} of $G\times \Sigma_n$ is a subgroup $\Gamma \subseteq G\times \Sigma_n$ that intersects $\{1_G\}\times \Sigma_n$ trivially.
\end{definition}

We recall the definition of an $N_{\infty}$-operad in $G$-simplicial sets.

\begin{definition}[\protect{\cite[Definition 2.6]{Rub1}}]
	An \emph{$N_{\infty}$-operad} is an operad in the category $\gsset$ such that: 
	\begin{enumerate}
		\item for every $n\ge 0$, the action of $\Sigma_n$ on $\ma{O}_n$ is free;
		\item for every graph subgroup $\Gamma \subseteq G\times \Sigma_n$, the simplicial set $(\ma{O}_n)^{\Gamma}$ is either contractible or empty;
		\item the simplicial sets $(\ma{O}_2)^G$ and $(\ma{O}_0)^G$ are non-empty.
	\end{enumerate}
\end{definition}

Recall now that graph subgroups can be characterised in terms of finite equivariant sets over subgroups of $G$.

\begin{lemma}[\protect{\cite[Lemma 2.4]{Rub1}}]
	For every graph subgroup $\Gamma\subseteq G\times \Sigma_n$, there is a unique subgroup $H\subseteq G$ and a homomorphism $\sigma:H\to \Sigma_n$ such that $\Gamma=\{(h,\sigma(h))~|~h\in H\}$.
	
\end{lemma}

For every finite $H$-set $T$, we choose once and for all an ordering $T\simeq \{1,\ldots, |T|\}$. The previous lemma gives a correspondence between finite $H$-sets and graph subgroups.
If $H\subseteq G$ is a subgroup of $G$ and $T$ a finite $H$-set, we denote by $\Gamma_T\subseteq G\times \Sigma_{|T|}$ the graph subgroup associated with the homomorphism $\sigma:H\to \Sigma_{|T|}$ that corresponds to $T$ through the chosen ordering on $T$.
If $\ma{O}$ is an $N_{\infty}$-operad, then we consider, for every subgroup $H\subseteq G$, the following collection
$$
\Ind_{\ma{O}} (H)=\left\{T\in \hset{H}~|~\left(\ma{O}_{|T|}\right)^{\Gamma_T}\neq \emptyset \right\}.
$$
As proven in \cite[Theorem 1.2]{BH1}, these collections of equivariant finite sets form a $G$-indexing system in the sense of the following definition.

\begin{definition}[\protect{\cite[Definition 2.12]{Rub1}}]\label{ind} 
	A $G$-\emph{indexing system} $\mathcal{I}$ consists, for every subgroup $H\subseteq G$, of a class $\mathcal{I}(H)$ of finite $H$-sets satisfying the following conditions:
	\begin{enumerate}
		\item for every subgroup $H\subseteq G$, the class $\mathcal{I}(H)$ contains all trivial $H$-sets;
		\item for every subgroup $H\subseteq G$ and finite $H$-sets $A$ and $B$ such that $B$ is isomorphic to $A$ in $\hset{H}$, if $B\in \mathcal{I}(H)$, then $A\in \mathcal{I}(H)$;
		\item for every sequence of subgroups $K\subseteq H\subseteq G$ and finite $H$-set $A$, if $A\in \mathcal{I}(H)$, then $\res{H}{K}A\in \mathcal{I}(K)$;
		\item for every subgroup $H \subseteq G$, $g\in G$ and finite $H$-set $A$, if $A\in \mathcal{I}(H)$, then $\operatorname{c}_gA\in \mathcal{I}(H^g)$;
		\item for every subgroup $H\subseteq G$ and finite $H$-sets $A$ and $B$, if $B$ is an $H$-equivariant subset of $A$ and $A\in \mathcal{I}(H)$, then $B\in \mathcal{I}(H)$;
		\item for every subgroup $H\subseteq G$ and finite $H$-sets $A$ and $B$, if $A\in \mathcal{I}(H)$ and $B\in \mathcal{I}(H)$, then $A\sqcup B \in \mathcal{I}(H)$;
		\item for every sequence of subgroups $K\subseteq H \subseteq G$ and finite $K$-set $A$, if $A\in \mathcal{I}(K)$ and $H/K \in \mathcal{I}(H)$, then $\ind{H}{K}A\in \mathcal{I}(H)$.
	\end{enumerate}
	Given a $G$-indexing system $\Ind$, we say that a finite $H$-set $T$ is \emph{admissible} if it belongs to $\Ind (H)$. The $G$-indexing systems form a poset under the inclusion that we denote by $\IndG$. 
\end{definition}

The $N_{\infty}$-operads for $G$ form an $\infty$-category that can be defined as follows.

\begin{definition}[\protect{\cite[Definition 3.9]{BH1}}]
	A morphism $f:\ma{O} \to \ma{P}$ of operads in $\gsset$ is a \emph{graph weak equivalence of $G$-operads} if for every graph subgroup $\Gamma \subseteq G\times \Sigma_n$, the morphism $f^{\Gamma}:{\ma{O}_n}^{\Gamma} \to {\ma{P}_n}^{\Gamma}$ is a weak homotopy equivalence of simplicial sets. We denote by $\gOpe$ the $\infty$-category obtained from $\gope$ by inverting graph weak equivalences and respectively by $\nope$ and $\Nope$ the full subcategories of $\gope$ and $\gOpe$ spanned  by $N_{\infty}$-operads.
\end{definition}

The $\infty$-category of $N_{\infty}$-operads is equivalent to the poset of $G$-indexing systems.

\begin{theorem}
	The functor $\nope \to \IndG$ that sends an $N_{\infty}$-operad $\ma{O}$ to its associated $G$-indexing system $\Ind_{\ma{O}}$ sends graph weakly equivalent $N_{\infty}$-operads to the same $G$-indexing system. Moreover, the induced functor $\Nope \to \IndG$ is an equivalence of $\infty$-categories. 
	\begin{proof}
		Combining \cite[Theorem 1.2]{BH1} and \cite[Proposition 5.5]{BH1}, one deduces that the functor \mbox{$\nope \to \IndG$} sends graph weakly equivalent $N_{\infty}$-operads to the same $G$-indexing system and that the induced functor $\Nope \to \IndG$ is fully faithful. The fact that this functor is essentially surjective follows independently from \cite[Theorem 7.2]{Rub1} and \cite[Corollary IV]{BonPer1}.
	\end{proof}
\end{theorem}

\subsection{The \texorpdfstring{$\infty$-categories}{infinity categories} of \texorpdfstring{$\cO$-algebras}{O-algebras} for an \texorpdfstring{$N_{\infty}$}{N}-operad \texorpdfstring{$\cO$}{O}}
Let $\cO$ be an $N_{\infty}$-operad. In this subsection, we define the $\infty$-category $\textbf{Alg}_{\cO}$ of $\cO$-algebras. We recall first the definition of $G$-weak homotopy equivalences and of the $\infty$-category of $G$-equivariant spaces.

\begin{definition}
	A morphism $f:X\to Y$ in $\gsset$ is a \emph{$G$-weak homotopy equivalence} if, for every subgroup $H\subseteq G$, the morphism $f^H:X^H\to Y^H$ is a weak homotopy equivalence in $\sset$. The \emph{$\infty$-category of $G$-equivariant spaces} $\gs$ is the $\infty$-category obtained from $\gsset$ by inverting $G$-weak homotopy equivalences. 
\end{definition}
We now define the $\infty$-category of $\ma{O}$-algebras.
\begin{definition}
A morphism $f:X\to Y$ in $\Alg_{\ma{O}}\left(\gsset\right)$ is a \emph{weak equivalence of $\ma{O}$-algebras} if $U(f):U(X)\to U(Y)$ is a $G$-weak homotopy equivalence with $U:\Alg_{\ma{O}}\left(\gsset\right)\to \gsset$ the forgetful functor. 
The \emph{$\infty$-category of $\ma{O}$-algebras}  $\mathbf{Alg}_{\ma{O}}$ is the $\infty$-category obtained from $\Alg_{\ma{O}}\left(\gsset\right)$ by inverting weak equivalences of $\ma{O}$-algebras.
\end{definition}

The $\infty$-category $\mathbf{Alg}_{\ma{O}}$ admits the following useful description.

\begin{proposition}
	The category $\Alg_{\ma{O}}\left(\gsset\right)$ admits a combinatorial and simplicial model structure right-induced from the genuine model structure on $\gsset$ through the free-forgetful adjunction
	$$
	F:\gsset \rightleftarrows \Alg_{\ma{O}}\left(\gsset\right):U.
	$$
	Moreover, the $\infty$-category $\mathbf{Alg}_{\ma{O}}$ is equivalent to the homotopy coherent nerve of the simplicial subcategory of $\Alg_{\ma{O}}\left(\gsset\right)$ spanned by objects that are fibrant and cofibrant. 
	\begin{proof}
		The category $\gsset$ is a monoidal model category with a cofibrant unit and with a symmetric monoidal fibrant replacement functor given by $\operatorname{Ex}^{\infty}$. Moreover, the monoidal structure is given by the cartesian product and it follows that the classical interval object $\Delta^1$ is endowed with a unique structure of cocommutative coalgebra. It follows from \cite[Theorem 2.1]{MR2342815} that the right-induced model structure on $\Alg_{\ma{O}}\left(\gsset\right)$ exists. If $X$ is an $\ma{O}$-algebra and if $K$ is a simplicial set, then the power $X^{K}$ in $\gsset$ is naturally endowed with the structure of an $\ma{O}$-algebra and this construction induces a structure of simplicial model category on $\Alg_{\ma{O}}\left(\gsset\right)$. Finally, the second part of the proposition follows from \cite[Theorem 1.3.4.20]{HA}.
	\end{proof}
\end{proposition}

\begin{remark}\label{topegsimp}
	In this thesis, we chose to work with $G$-simplicial sets instead of $G$-topological spaces. There is a monoidal Quillen equivalence
	$$
	\left|-\right|:\gsset \rightleftarrows \gtop:\Sing
	$$
	between $\gsset$ and $\gtop$ respectively, endowed with the genuine model structures.
    The functor $\left|-\right|$ preserves finite limits and it follows that the operad $\left|\ma{O}\right|$ is an $N_{\infty}$-operad in $\gtop$ with the same $G$-indexing system as $\cO$. Moreover, the functor $\left|-\right|:\gsset \to \gtop$ preserves weak equivalences and it follows from \cite[Theorem 4.7]{BM03} that the $\infty$-categories $\mathbf{Alg}_{\ma{O}}$ and $\mathbf{Alg}_{\left|\ma{O}\right|}$ are equivalent. 
\end{remark}

Let $\Ind$ be a $G$-indexing system. Different $N_{\infty}$-operads that correspond to $\Ind$ should produce equivalent $\infty$-categories, and one can make the following conjecture.

\begin{conjecture}\label{equininf}
	Let $\ma{O}$ and $\ma{P}$ be two $N_{\infty}$-operads associated with $\Ind$. The $\infty$-categories $\mathbf{Alg}_{\ma{O}}$ and $\mathbf{Alg}_{\ma{P}}$ are equivalent.
\end{conjecture}

Note however that $N_{\infty}$-operads are not $\Sigma$-cofibrant, as pointed out in \cite[Proposition 2.1]{greg3}. In particular, many results of \cite{BM03} do not apply to $N_{\infty}$-operads and, in particular, the proof of \cite[Theorem 4.4]{BM03} does not directly generalise. In this paper, we will work with Rubin's combinatorial $N_{\infty}$-operad $\EII$ as a model, and we will recall its definition in \ref{rubmod}. For the rest of the paper, we will regard $\mathbf{Alg}_{\EII}$ as the \emph{$\infty$-category of $N_{\infty}$-algebras for $\Ind$}. Note however that, as long as \ref{equininf} is not proven, the results of this paper only apply to algebras over $\EII$.

\section{Normed symmetric monoidal categories}\label{Norcat}
In this section, we recall Rubin's notion of $\Ind$-normed symmetric monoidal categories, and we give an explicit description of the operad $\EI$ in $\catG$ whose algebras in $\catG$ are $\Ind$-normed symmetric monoidal categories. The operad $\EI$ lives in the category of $G$-groupoids and its nerve $\EII:=N(\EI)$ is an $N_{\infty}$-operad with associated $G$-indexing system $\Ind$. 
The description that we give of $\EI$ is equivalent, but slightly different from the one given by Rubin in \cite{Rub1} and \cite{Rub2}. Our description will be useful to give an explicit construction of the free $\EII$-algebras over finite $G$-sets in terms of $\Ind$-normed symmetric monoidal groupoids, which is a part of the construction of the equivalence of \ref{THETH}. 

\subsection{Preliminaries on \texorpdfstring{$\Ind$}{I}-normed symmetric monoidal categories}
In this subsection, we recall Rubin's notion of $\Ind$-normed symmetric monoidal categories. We refer the reader to \cite{Rub2} and \cite{Rub3} for more details on the topic.
\begin{definition}
	Let $H\subseteq G$ be a subgroup, $T$ a finite $H$-set and $\cC$ a $G$-category. Denote by $\cC^{\times T}$ the $H$-category given by the $|T|$-fold cartesian power $\cC^{|T|}$ endowed with the $H$-action given on objects by
	$$
	h\cdot \left(X_1,\ldots,X_{|T|}\right):=\left(h X_{\sigma(h)^{-1}1},\ldots ,h X_{\sigma(h)^{-1}|T|}\right)
	$$
	and by the same formula on morphisms, where $\sigma:H\to \Sigma_{|T|}$ is the homomorphism that corresponds to the action on $T$ through the bijection $T\simeq \{1,\ldots, |T|\}$ given by the chosen ordering on $T$.
	A \emph{$T$-external norm} is an $H$-functor $\cC^{\times T}\to \cC$.
\end{definition}

\begin{definition}[\protect{\cite[Definition 2.3]{Rub2}}]
	An \emph{$\Ind$-normed symmetric monoidal category} $\cC$ is a $G$-symmetric monoidal category $(\cC,\otimes,e,\alpha,\lambda,\rho,\beta)$ together with:
	\begin{enumerate}
		\item a $T$-external norm $\bigotimes_{T}:\cC^{\times T}\to \cC$ for every subgroup $H\subseteq G$ and every admissible $H$-set $T$;
		\item (\emph{untwistors}) a nonequivariant natural isomorphism 
		$$
		\textstyle v_T:\bigotimes\nolimits_T(X_1,\ldots,X_{|T|})\to \bigotimes\nolimits_{|T|} \left(X_1,\ldots,X_{|T|}\right)
		$$
		for every admissible $H$-set $T$ such that for every $h\in H$, the diagram 
		$$
		\begin{tikzcd}
		h \bigotimes_T\left(X_1,\ldots,X_{|T|}\right) \arrow[r,"\Id{}"] \arrow[dd,"h v_T"'] & \bigotimes_T\left(h X_{\sigma(h)^{-1}1},\ldots,h  X_{\sigma(h)^{-1}|T|}\right) \arrow[d,"v_T"]\\
		\text{} & \bigotimes_{|T|}\left(h X_{\sigma(h)^{-1}1},\ldots,h  X_{\sigma(h)^{-1}|T|}\right) \arrow[d,"\sigma(h)^{-1}"] \\
		h \bigotimes_{|T|}\left(X_1,\ldots,X_{|T|}\right) \arrow[r,"\Id{}"'] & \bigotimes_{|T|}\left(hX_1,\ldots,h X_{|T|}\right)
		\end{tikzcd}
		$$
		commutes where $\bigotimes_{|T|} \left(X_1,\ldots,X_{|T|}\right)$ denotes the $|T|$-fold tensor product and where the morphism $\sigma(h)^{-1}$ denotes the canonical isomorphism given by the symmetric monoidal structure of $\cC$ that permutes the factors of $\bigotimes_{|T|}$ by $\sigma(h)^{-1}$.
	\end{enumerate}
\end{definition}

\begin{definition}[\protect{\cite[Definition 2.5]{Rub2}}]
	Let $\cC$ and $\cD$ be two $\Ind$-normed symmetric monoidal categories. A \emph{lax $\Ind$-normed functor} $\phi:\cC \to \cD$ consists of the following data:
	\begin{enumerate}
		\item a $G$-functor $\phi:\cC \to \cD$;
		\item a $G$-fixed morphism $\phi_e:e^{\cD} \to \phi e^{\cC}$;
		\item a $G$-natural transformation $\phi_{\otimes}:\phi X\otimes^{\cD}\phi Y\to \phi \left(X\otimes^{\cC} Y\right)$;
		\item for every admissible $H$-set $T$, an $H$-natural transformation
		\[ 
		\textstyle \phi_{\bigotimes_T}:\bigotimes\nolimits_T^{\cD}\left(\phi X_1,\ldots,\phi X_{|T|}\right) \to \phi \left(\bigotimes\nolimits_T^{\cC}\left(X_1,\ldots,X_{|T|}\right)\right)
		\]
		such that the usual lax symmetric monoidal diagrams relating $\alpha, \lambda, \rho, \beta$ to the maps 
		$\phi_e$ and $\phi_{\bigotimes}$ commute and such that the square
		$$
		\begin{tikzcd}[column sep=huge,row sep=large]
		\bigotimes_T^{\cD}\left(\phi X_1,\ldots, \phi X_{|T|}\right) \arrow[r,"v^{\cD}_T"] \arrow[d,"\phi _{\bigotimes_T}"']  & \bigotimes_{|T|}^{\cD}\left(\phi X_1,\ldots,\phi X_{|T|}\right) \arrow[d,"\phi_{\bigotimes_{|T|}}"] \\
		\phi \left(\bigotimes_T^{\cC}\left(X_1,\ldots,X_{|T|}\right)\right) \arrow[r,"\phi(v_T^{\cC})"'] & \phi \left(\bigotimes_{|T|}^{\cC}\left(X_1,\ldots,X_{|T|}\right)\right)
		\end{tikzcd}
		$$
		commutes for every admissible $H$-set $T$ where $\phi_{\bigotimes_{|T|}}$ is the map constructed with iterations of $\phi_{\bigotimes}$ and $\phi_e$. A lax $\Ind$-normed functor is a \emph{strong $\Ind$-normed functor} (respectively \emph{strict $\Ind$-normed functor}) if the natural transformations $\phi_e$, $\phi_{\bigotimes}$ and $\phi_{\bigotimes_T}$ are isomorphisms (respectively identities).
	\end{enumerate}
\end{definition}

If $\phi:\cC \to \cD$ and $\psi:\cD \to \cE$ are lax $\Ind$-normed functors, then their composition $\psi \circ \phi$ is also canonically a lax $\Ind$-normed functor. The morphism $(\psi \circ \phi)_e$ is given by the composition 
$$
e_{\cE} \overset{\psi_e}{\longrightarrow} \psi e_{\cD} \overset{\psi (\phi_e)}{\longrightarrow} (\psi\circ \phi) e_{\cC},
$$
the morphism $(\psi\circ \phi)_{\otimes}$ by the composition
$$
((\psi \circ \phi) X)\otimes^{\cE} ((\psi \circ \phi) Y) \overset{\psi_{\otimes}}{\longrightarrow} \psi \left(\phi X \otimes^{\cD} \phi Y\right) \overset{\psi (\phi_{\otimes})}{\longrightarrow} (\psi \circ \phi) \left(X \otimes^{\cC} Y\right)
$$
and, for every admissible finite $H$-set $T$, the morphism $(\psi\circ \phi)_{\otimes_T}$ is given by the composition of the morphism
$$
\textstyle\bigotimes_T\left((\psi \circ \phi) X_1 ,\ldots, (\psi \circ \phi) X_{|T|}\right) \overset{\psi_{\bigotimes_T}}{\longrightarrow} \psi \left(\bigotimes_T \left(\phi X_1 ,\ldots , \psi X_{|T|}\right)\right) 
$$
with the morphism
$$
\textstyle \psi \left(\bigotimes_T \left(\phi X_1 ,\ldots , \psi X_{|T|}\right)\right) \overset{\psi (\phi_{\bigotimes_T})}{\longrightarrow} (\psi \circ \phi) \left(\bigotimes_T \left(X_1,\ldots, X_{|T|}\right)\right).
$$
If $\phi$ and $\psi$ are strong $\Ind$-normed functors (resp. strict), then $\psi \circ \phi$ is also a strong $\Ind$-normed functor (resp. strict).
\begin{definition}[\protect{\cite[Definition 2.6]{Rub2}}]
	Let $\cC$ and $\cD$ be $\Ind$-normed symmetric monoidal categories and $\phi,\psi:\cC \rightrightarrows \cD$ a pair of lax $\Ind$-normed functors. An \emph{$\Ind$-normed natural transformation} $\omega:\phi \Rightarrow \psi$ is a $G$-natural transformation $\omega:\phi \Rightarrow \psi$ such that the usual monoidal transformation squares relating $\phi_e$, $\psi_e$, $\phi_{\bigotimes}$ and $\psi_{\bigotimes}$ to $\omega$ commute and such that the square 
	$$
	\begin{tikzcd}[column sep=huge,row sep=large]
	\bigotimes_T^{\cD}\left(\phi X_1,\ldots ,\phi X_{|T|}\right) \arrow[r,"\bigotimes_T(\omega\text{,}\ldots\text{,}\omega)"] \arrow[d,"\phi_{\bigotimes_T}"'] & \bigotimes_T^{\cD}\left(\psi X_1,\ldots, \psi X_{|T|}\right)
	\arrow[d,"\psi_{\bigotimes_T}"] \\ \phi \left(\bigotimes_T^{\cC}\left(X_1,\ldots,X_{|T|}\right)\right) \arrow[r,"\omega"'] & \psi\left(\bigotimes_T^{\cC}\left(X_1,\ldots,X_{|T|}\right)\right)
	\end{tikzcd}
	$$
	commutes for every admissible $H$-set $T$.
\end{definition}

\begin{notation}
	The structure of a 2-category on $\catG$ lifts to $\Ind$-normed symmetric monoidal categories and we denote by $\NSSlaxI$ (resp. $\NSStMI$ and $\NSSMI$) the $2$-category of $\Ind$-normed symmetric monoidal categories, lax $\Ind$-normed functors (resp. strong and strict)
	and $\Ind$-normed natural transformations.
\end{notation}

The category of strong $\Ind$-normed functors $\NSStMI (\cC,\cD)$ between two $\Ind$-normed symmetric monoidal categories $\cC$ and $\cD$ is naturally endowed with a structure of a symmetric monoidal category.

\begin{proposition}\label{2sem}
	Let $\phi,\psi:\cC \rightrightarrows \cD$ be two strong $\Ind$-normed functors. The composition
	\[
	\phi \otimes \psi : \cC \overset{\Delta}{\longrightarrow} \cC \times \cC \overset{\phi\times \psi}{\longrightarrow} \cD \times \cD \overset{\otimes_{\cD}}{\longrightarrow} \cD
	\]
	is naturally endowed with the structure of a strong $\Ind$-normed functor. Moreover, if $\tau_1:\phi_1 \Rightarrow \phi_2$ and $\tau_2:\psi_1 \Rightarrow \psi_2$ are two $\Ind$-normed natural transformations between strong $\Ind$-normed functors, the natural transformation $\tau_1 \otimes \tau_2$ defined as $\bigotimes_{\cD}\circ (\tau_1\times \tau_2)\circ \Delta$ is an $\Ind$-normed natural transformation between $\phi_1 \otimes \psi_1$ and $\phi_2 \otimes \psi_2$. Given two $\Ind$-normed symmetric monoidal categories $\cC$ and $\cD$, these constructions endow $\NSStMI (\cC,\cD)$ with the structure of a symmetric monoidal category. 
	\begin{proof}
		The $G$-fixed isomorphism $(\phi \otimes \psi)_e$ is the composition 
		\[
		e_\cD \simeq e_\cD \otimes e_\cD\overset{\phi_e\otimes \psi_e}{\longrightarrow}\phi (e_\cC) \otimes \psi (e_\cC),
		\]
		the $G$-natural transformation
		$(\phi\otimes \psi)_{\otimes}$ is given by the composition of the morphism
		\[
		(\phi (X) \otimes \psi (X))\otimes (\phi (Y) \otimes \psi (Y)) \simeq (\phi (X)\otimes \phi (Y) ) \otimes (\psi (X)\otimes \psi (Y)) 
		\]
        with the morphism
        $$
        (\phi (X)\otimes \phi (Y) ) \otimes (\psi (X)\otimes \psi (Y)) \overset{(\phi_\otimes)\otimes (\psi_{\otimes})}{\simeq}
		\phi(X\otimes Y)\otimes \psi(X\otimes Y)
        $$
		and, for every admissible $H$-set $T$, the $H$-natural transformation $(\phi \otimes \psi )_{\bigotimes_T}$ is given by the composition 
		\begin{align*} 
		&\textstyle\bigotimes_T \left(\phi\left(X_1\right)\otimes \psi\left(X_1\right), \ldots, \phi \left(X_{|T|}\right)\otimes \psi\left(X_{|T|}\right)\right) \\
		\overset{\alpha}{\simeq} &\textstyle \left(\bigotimes_T \left(\phi \left(X_1\right),\ldots, \phi\left(X_{|T|}\right)\right) \right)\otimes \left(\bigotimes_T \left(\psi \left(X_1\right),\ldots, \psi\left(X_{|T|}\right)\right)\right) \\
		\overset{\left(\phi_{\bigotimes_T}\right)\otimes \left(\psi_{\bigotimes_T}\right)}{\simeq} &\textstyle \phi \left(\bigotimes_T \left(X_1,\ldots, X_{|T|}\right)\right)\otimes \psi \left(\bigotimes_T \left(X_1,\ldots, X_{|T|}\right)\right)
		\end{align*}
		where $\alpha$ is the $H$-natural transformation given by the composition 
		\begin{align*}
		&\textstyle \bigotimes_T \left(\phi\left(X_1\right)\otimes \psi\left(X_1\right), \ldots, \phi \left(X_{|T|}\right)\otimes \psi\left(X_{|T|}\right)\right) \\
		\overset{v_T}{\simeq} &\textstyle \bigotimes_{|T|} \left(\phi\left(X_1\right)\otimes \psi \left(X_1\right), \ldots, \phi \left(X_{|T|}\right)\otimes \psi\left(X_{|T|}\right)\right) \\
		\simeq &\textstyle \left(\bigotimes_{|T|} \left(\phi \left(X_1\right),\ldots, \phi \left(X_{|T|}\right)\right) \right)\otimes \left(\bigotimes_{|T|} \left(\psi \left(X_1\right),\ldots, \psi \left(X_{|T|}\right)\right)\right)  \\
		\overset{{v_T}^{-1} \otimes {v_T}^{-1}}{\simeq} &\textstyle \left(\bigotimes_T \left(\phi \left(X_1\right),\ldots, \phi\left(X_{|T|}\right)\right) \right)\otimes \left(\bigotimes_T \left(\psi \left(X_1\right),\ldots, \psi\left(X_{|T|}\right)\right)\right).
		\end{align*}
		A direct inspection shows that $\phi\otimes \psi$ is a strong $\Ind$-normed functor and that the result also stands on $\Ind$-normed natural transformations. 
		The structural morphisms of the symmetric monoidal category $\NSStMI(\cC,\cD)$ are given componentwise by the ones of $\cD$.
	\end{proof}
\end{proposition}

\begin{remark}
	Note that if $\phi,\psi:\cC \rightrightarrows \cD$ is a pair of strict $\Ind$-normed functors, then the strong $\Ind$-normed functor $\phi \otimes \psi$ is not strict in general. In particular, the previous proposition fails for strict $\Ind$-normed functors.
\end{remark}

It follows from \ref{2sem} that the $2$-category $\NSStMI$ is a $\mathbf{SMC}$-category in the sense of \cite[Definition 2.1]{MR2770075}. Moreover, $\NSStMI$ admits finite products and is, in particular, a semi-additive $2$-category. 

\begin{proposition}\label{2sa}
	The $2$-category $\NSStMI$ is semi-additive, ie, $\NSStMI$ admits direct sums.
	\begin{proof}
		Let $\cC$ and $\cD$ be two $\Ind$-normed symmetric monoidal categories. The product $\cC \times \cD$ is canonically endowed with a structure of $\Ind$-normed symmetric monoidal category which is the cartesian product of $\cC$ and $\cD$ in the $2$-category $\NSStMI$. 
		The final $G$-category $*$ is also canonically an $\Ind$-normed symmetric monoidal category and is a zero object in the $2$-category $\NSStMI$. If we denote by $0_{\cC,\cD}:\cC \to \cD$ the strong $\Ind$-normed functor constant on $e_{\cD}$,
		then we have to show that $\cC \times \cD$ endowed with the strong $\Ind$-normed functors $i_{\cC}=(\Id{\cC},0_{\cC,\cD}):\cC \to \cC \times \cD$ and $i_{\cD}=(0_{\cD,\cC},\Id{\cD}):\cD \to \cC\times \cD$ is the coproduct of $\cC$ and $\cD$ in $\NSStMI$.
		More precisely, if we denote by $p_{\cC}:\cC\times \cD \to \cC$ and $p_{\cD}:\cC \times \cD \to \cD$ the projections, then we prove that the functor
		$$
		\left((-)\circ p_{\cC}\right)\otimes ((-)\circ p_{\cD}):\NSStMI (\cC,\cE)\times \NSStMI (\cD,\cE) \to \NSStMI (\cC\times \cD,\cE)
		$$
		that sends a pair of strong $\Ind$-normed functors $\phi:\cC \to \cE$ and $\psi:\cD \to \cE$ to $(\phi \circ p_{\cC})\otimes (\psi \circ p_{\cD})$ is an inverse of the functor 
		$$
		((-)\circ i_{\cC},(-) \circ i_{\cD}):\NSStMI (\cC \times \cD, \cE) \to \NSStMI (\cC, \cE) \times \NSStMI (\cD, \cE).
		$$
		Given a pair of strong $\Ind$-normed functors  $\phi:\cC \to \cE$ and $\psi:\cD \to \cE$, we have $((\phi \circ p_{\cC})\otimes (\psi \circ p_{\cD}))\circ i_{\cC} =\phi \otimes 0_{\cC,\cE} \simeq \phi$ and $((\phi \circ p_{\cC})\otimes (\psi \circ p_{\cD}))\circ i_{\cD} = 0_{\cD,\cE} \otimes \psi  \simeq \psi$ and if $\zeta :\cC \times \cD \to \cE$ is a strong $\Ind$-normed functor, then 
		$(\zeta \circ i_{\cC} \circ p_{\cC})\otimes (\zeta \circ i_{\cD} \circ p_{\cD})$ is the strong $\Ind$-normed functor that sends an object $(X,Y)$ of $\cC \times \cD$ to $\zeta (X,e_{\cD})\otimes \zeta(e_{\cC},Y) \overset{\zeta_{\otimes}}{\simeq} \zeta\left(X\otimes e_{\cC},Y\otimes e_{\cD}\right)\simeq \zeta(X,Y)$ and the result follows. 
	\end{proof}
\end{proposition}

\subsection{Rubin's combinatorial models of \texorpdfstring{$N_{\infty}$-operads}{N-infinity operads}}\label{rubmod}

In this subsection, we give an explicit description of the operad $\EI$ in $\catG$ defined in \cite{Rub1} and \cite{Rub2} whose algebras in $\catG$ are $\Ind$-normed symmetric monoidal categories. 
We then use this description to give a short proof of the fact that the operad $\EII=N(\EI)$ is an $N_{\infty}$-operad associated with our chosen $G$-indexing system $\Ind$. 

\begin{construction}\label{transg}
	If $X$ is a $G$-set, then the \emph{translation $G$-groupoid $\widetilde{X}$ of $X$} is the $G$-category defined by the following diagram 
	$$
	\begin{tikzcd}
	(X\times X) \times_X (X\times X)\simeq X\times X \times X \arrow[r,"(p_1\text{, }p_3)"]& X\times X \arrow[r,shift right=0.25cm,"p_2"'] \arrow[r,shift left=0.3cm,"p_1"] &X. \arrow[l,"\Delta"']
	\end{tikzcd}
	$$
	The associated functor $\widetilde{(-)}:\gsett \to \catG$ is a right adjoint of the object functor $\Ob:\catG \to \gsett$ and in particular, preserves limits. It follows that if $\ma{O}$ is an operad in $\gsett$, then $\widetilde{\ma{O}}$ is canonically an operad in $\catG$.
\end{construction}

\begin{definition}
	The operad $\EI$ is the operad $\widetilde{\SN}$ with $\SN$ the free operad in the category $\gsett$ over the symmetric sequence $S_{\Ind}$ defined by
	$$
	(S_{\Ind})_n=\bigsqcup_{T\in \Ind \text{, } |T|=n} \frac{G\times \Sigma_n}{\Gamma
		_T}.
	$$
\end{definition}

We will now give an explicit description of $\SN$. The construction of the free operad over a symmetric sequence in $\gsett$ can be made explicit in general, as in \cite[Section 5]{Rub1} for instance, though its general description can be complicated to handle. 
However, if we denote by $V_{\Ind}$ the non-symmetric sequence of $G$-sets defined by 
$$
(V_{\Ind})_n=\bigsqcup_{\substack{H\subseteq G,~ |T|=n, \\ T\in \Ind(H)}}G/H ,
$$
the choice for every subgroup $H\subseteq G$ of a set of coset representatives of $G/H$ induces a nonequivariant isomorphism of symmetric sequences $S_{\Ind}\simeq V_{\Ind}\times \Sigma_*$ where $\Sigma_*$ denotes the symmetric sequence of sets given by the symmetric groups. If we forget the action of $G$, then it follows that the operad $\SN$ can be described as the product $\FVI \times \Sigma_*$ with $\FVI$ the free nonsymmetric operad over the sequence $V_{\Ind}$. 
We recall now the construction of the nonsymmetric operad $\FVI$ in $\set$ in terms of planar rooted trees and we then give an explicit description of the required $G$-action on $\SN$. 

\begin{remark}
	The freeness of the symmetric sequence of sets $S_{\Ind}$ was already used in \cite[Theorem 5.21]{Rub1} to give an explicit description of $F(S_{\Ind})$. However, \cite[Theorem 5.21]{Rub1} does not construct $F(S_{\Ind})$ as the product $\FVI \times \Sigma_*$. We use this presentation to construct explicitly the free $N_{\infty}$-algebras over finite $G$-sets in \ref{exa}.
\end{remark}

\begin{construction}
	For every subgroup $H\subseteq G$, choose a set of coset representatives $\{1_G=g_1,\ldots, g_{[G:H]}\}$ of $G/H$. Consider the formal symbols $r\bigotimes_T$ and $\Id{}$ with $r$ a $G/H$ coset representative and $T$ an admissible $H$-set. A $\emph{formal external norm}$ is a formal composite of the previous symbols defined by recursion as follows:
	\begin{enumerate}
		\item the symbol $\Id{}$ is a formal external norm;
		\item if $\theta_1,\ldots, \theta_{|T|}$ are formal external norms, then $r\bigotimes_T \left(\theta_1,\ldots, \theta_{|T|}\right)$ is a formal external norm.
	\end{enumerate}
	If $\theta$ is a formal external norm, then its \emph{length} $|\theta|$ is defined by the following recursion:
	\begin{enumerate}
		\item if $\theta=\Id{}$, then $|\theta|=1$;
		\item if $\theta=r\bigotimes_T\left(\theta_1,\ldots, \theta_{|T|}\right)$, then $|\theta|=|\theta_1|+\cdots+ |\theta_{|T|}|$.
	\end{enumerate}
	
	The nonsymmetric operad $\FVI$ in $\set$ can now be described as the operad whose set of n-ary operations $\FVI_n$ is the one of formal external norms of length $n$. The identity of $\FVI$ is $\Id{}$ and if $\theta$ and $\tau_1, \ldots, \tau_n$ are formal external norms with $|\theta|=n$, their composition $\gamma \left(\theta,\tau_1,\ldots,\tau_n\right)$ is defined by recursion as follows:
	\begin{enumerate}
		\item if $\theta=\Id{}$, then $\gamma (\theta,\tau_1)=\tau_1$;
		\item if $\theta=r\bigotimes_T(\theta_1,\ldots, \theta_k)$, then $$\textstyle \gamma (\theta,\tau_1,\ldots,\tau_n)=r\bigotimes_T\left(\gamma\left(\theta_1,\tau_1,\ldots,\tau_{|\theta_1|}\right),\ldots, \gamma\left(\theta_k,\tau_{|\theta|-|\theta_k|},\ldots,\tau_n\right)\right).$$
	\end{enumerate}
	
	There is an action of $G$ on formal external norms defined by recursion as follows:
	\begin{enumerate}
		\item if $\theta=\Id{}$, then $g\cdot \theta=\theta$;
		\item if $\theta=r\bigotimes_T\left(\theta_1,\ldots, \theta_{|T|}\right)$, then $g\cdot \theta =r'\bigotimes_T \left(h\cdot \theta_{\sigma(h)^{-1}1},\ldots, h\cdot\theta_{\sigma(h)^{-1}|T|}\right)$ with $r'$ the chosen coset representative of $gr$ and $h$ the unique element of $H$ such that $gr=r'h$ and where $\sigma:H\to \Sigma_{|T|}$ is the homomorphism that corresponds to $T$ through the chosen ordering $T\simeq \{1,\ldots, |T|\}$.
	\end{enumerate}
	The $G$-action on formal external norms is compatible with the operadic structure of $\FVI$ and it follows that $\FVI$ is canonically an operad in $\gsett$.
\end{construction}

\begin{remark}
	As we already mentioned, $\FVI$ is the free nonsymmetric operad in $\set$ over the sequence 
	$$
	(V_{\Ind})_n=\bigsqcup_{\substack{H\subseteq G,~ |T|=n, \\ T\in \Ind(H)}}G/H ,
	$$
	and it follows that a formal external norm can be seen as a planar rooted tree with $n$-ary nodes decorated with the symbols $r\bigotimes_T$ with $T$ an admissible $H$-set of cardinality $n$, $r$ a chosen coset representative of $G/H$, where leaves represent occurrences of the identity $\Id{}$.
	The figure below shows an example
	
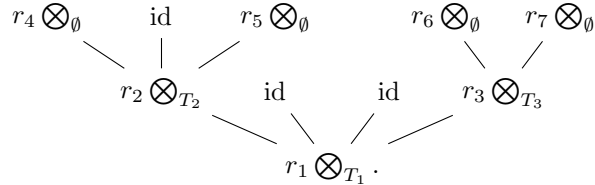
\begin{center}
		\begin{tikzpicture}[grow via three points={%
			one child at (0,1) and two children at (-0.75,1) and (0.75,1)}]
		\node at (0,-4.5) {$r_1\bigotimes_{T_1}.$} child{node {$r_2\bigotimes_{T_2}$}child{node{$r_4\bigotimes_{\emptyset}$} } child{node {$\Id{}$}} child{node {$r_5 \bigotimes_{\emptyset}$}}} child {node {$\Id{}$}} child{node {$\Id{}$}} child{node {$r_3\bigotimes_{T_3}$} child{ node {$r_6 \bigotimes_{\emptyset}$}} child{ node {$r_7 \bigotimes_{\emptyset}$}}} ;
		\end{tikzpicture}
	\end{center}
	Note that for the construction of $\FVI$, the choices for every subgroup $H$ of $G$, of orderings on admissible $H$-sets, and of coset representatives of $G/H$, play only a role in the definition of the $G$-action. The choices of coset representatives of $G/H$ correspond to the necessary choices to construct the nonequivariant isomorphism $S_{\Ind}\simeq V_{\Ind}\times \Sigma_*$ with $\Sigma_*$ the symmetric sequence associated with the symmetric groups. Moreover, note that $\FVI$ is the free operad over $V_{\Ind}$ only in the category $\set$ and not as an operad in $\gsett$.
\end{remark}

\begin{notation}
	If $\delta_1 \in \Sigma_{k_1}, \ldots, \delta_n\in \Sigma_{k_n}$ and $\delta\in \Sigma_n$ are permutations, then we denote by $\delta\num{\delta_1,\ldots,\delta_n}$ the permutation $\delta(k_1,\ldots,k_n) \circ \bigsqcup_{i=1}^n \delta_i \in \Sigma_{k_1 +\cdots + k_n}$ with $\delta(k_1,\ldots,k_n)$ the block permutation.
\end{notation}

We can now give the promised explicit description of the $G$-action on the operad $\SN$.
\begin{construction}
	Let $\theta$ be a formal external norm and $g$ an element of $G$. We define a permutation $\omega_{\theta}(g)$ in $\Sigma_{|\theta|}$ by recursion as follows:
	\begin{enumerate}
		\item if $\theta=\Id{}$, then $\omega_{\theta}(g)=\Id{\Sigma_{1}}$;
		\item if $\theta=r\bigotimes_{T}\left(\theta_1,\ldots,\theta_{|T|}\right)$, then $\omega_{\theta}(g)=\sigma(h)\num{\omega_{\theta_1}(h),\ldots,\omega_{\theta_{|T|}}(h)}$ with $h$ and $\sigma:H\to \Sigma_{|T|}$ given as in the definition of the $G$-action on formal external norms.
	\end{enumerate}
	The $G$-action on the operad $\SN$ is now defined on an element $(\theta,\delta)$ of $\SN_n=\FVI_n \times \Sigma_n$ by 
	$$
	g\cdot (\theta,\delta)= \left(g\cdot \theta, \omega_{\theta}(g)\circ \delta\right).
	$$
\end{construction}

\begin{remark}\label{actionwelldef}
	We can verify by induction that the action on formal external norms is well-defined. If $\theta$ is a formal external norm and $g$ and $g'$ are elements of $G$, we can check by induction that $g' (g \theta)=g'g  \theta$. First note that the result is trivial for $\theta=\Id{}$. Assume now that $\theta=r\bigotimes_T\left(\theta_1,\ldots,\theta_{|T|}\right)$ and that the result is true for $\theta_1, \ldots, \theta_{|T|}$. If we write $gr=r'h$ and $g'r'=r''h'$, we have
	\begin{align*} 
	g' (g \theta)=& \textstyle g' r'\bigotimes\nolimits_T \left(h \theta_{\sigma(h)^{-1}1},\ldots, h\theta_{\sigma(h)^{-1}|T|}\right) \\ 
	=&\textstyle r''\bigotimes\nolimits_T \left(h'h \theta_{\sigma(h'h)^{-1}1},\ldots, h'h\theta_{\sigma(h'h)^{-1}|T|}\right) \\
	=&g'g \theta. 
	\end{align*}
	
	A similar argument shows that $\omega_{g \theta}(g')\circ \omega_{\theta}(g)=\omega_{\theta}(g'g)$ and thus the action of $G$ on $\SN$ is also well defined. As a consequence, if $\theta$ is an $H$-fixed formal external norm, the restriction of $\omega_{\theta}:G\to \Sigma_{|\theta|}$ to $H$ is a homomorphism.
\end{remark}

\begin{remark}
	If $\theta$ is an $H$-fixed formal external norm, the homomorphism $\omega_{\theta}:H\to \Sigma_{|\theta|}$ depends on the chosen orderings on the admissible $H$-sets appearing in $\theta$ and on the choices of coset representatives of $G/H$. However, up to isomorphism, the corresponding $H$-set does not depend on these choices.
\end{remark}

\begin{notation}\label{equinorm}
	If $\theta$ is an $H$-fixed formal external norm, we denote by $T_{\theta}$ the finite $H$-set that corresponds to the homomorphism $\omega_{\theta}:H\to \Sigma_{|\theta|}$. 
\end{notation}

\begin{remark}
	An operation $(\theta,\delta)$ of $\SN$ can be seen as the planar rooted tree that corresponds to $\theta$ with an ordering on the leaves given by the permutation $\delta$. The picture below illustrates an example.
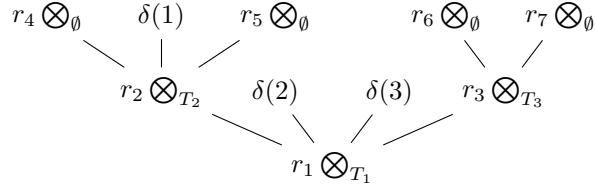
	\begin{center}
		\begin{tikzpicture}[grow via three points={%
			one child at (0,1) and two children at (-0.75,1) and (0.75,1)}]
		\node at (0,-4.5) {$r_1\bigotimes_{T_1}$} child{node {$r_2\bigotimes_{T_2}$}child{node{$r_4\bigotimes_{\emptyset}$} } child{node {$\delta(1)$}} child{node {$r_5\bigotimes_{\emptyset}$}}} child {node {$\delta(2)$}} child{node {$\delta(3)$}} child{node {$r_3\bigotimes_{T_3}$} child{ node {$r_6\bigotimes_{\emptyset}$}} child{ node {$r_7\bigotimes_{\emptyset}$}}} ;
		\end{tikzpicture}
	\end{center}
	If $\theta$ is a formal external norm and $g$ is an element of $G$, the permutation $\omega_{\theta}(g)$ corresponds to the permutation of the leaves obtained by acting on $\theta$ with $g$.
\end{remark}

\begin{notation}
	In what follows, if $T$ is an admissible finite $H$-set, then we denote by $\bigotimes_T$ the formal external norm $1_G \bigotimes_T(\Id{},\ldots, \Id{})$. If $\theta$ is any formal external norm, then we also denote by $\theta$ the operation $\left(\theta,\Id{\Sigma_{|\theta|}}\right)$ of $\SN$.
\end{notation}

In order to show that the description of the $G$-action of $\SN$ that we gave is the right one, we show now, under our description, that $\SN$ is the free operad in $\gsett$ over the symmetric sequence $S_{\Ind}$.

\begin{proposition}\label{freeop}
	The $G$-operad $\SN$ is the free operad in $\gsett$ over the symmetric sequence $S_{\Ind}$.
	\begin{proof}
		Note first that, given a symmetric sequence of $G$-sets $S$, a morphism $S_{\Ind}\to S$ is the data of, for every admissible $H$-set $T$, a $\Gamma_{T}$-fixed element of $S_{|T|}$.
		The operation $\bigotimes_T$ of $\SN$ is $\Gamma_T$-fixed and we thus get a morphism  $S_{\Ind}\to \SN$ of symmetric sequences of $G$-sets. Consider now any morphism of symmetric sequences $\phi:S_{\Ind} \to \mathcal{O}$ with $\mathcal{O}$ any $G$-operad. For every admissible $H$-set $T$, denote by $\bigotimes_T^{\mathcal{O}}$ the $\Gamma_T$-fixed operation of $\ma{O}$ that corresponds to $\phi$. If $\theta$ is a formal external norm, consider the operation $\Phi(\theta)$ of $\mathcal{O}$ defined by the following recursion:
		\begin{enumerate}
			\item if $\theta=\Id{}$, then $\Phi(\theta)=\Id{\mathcal{O}}$;
			\item if $\theta=r\bigotimes_T\left(\theta_1,\ldots,\theta_{|T|}\right)$, then $\Phi (\theta)=r \cdot  \gamma \left(\bigotimes_T^{\mathcal{O}},\Phi\left(\theta_1\right),\ldots,\Phi \left(\theta_{|T|}\right)\right)$.
		\end{enumerate}
		We can now define a morphism of $G$-operads $\Phi:\SN \to \mathcal{O}$ by setting $\Phi (\theta,\delta)= \Phi(\theta)\cdot \delta$.
		It is clear that $\Phi$ is a morphism of operads in $\set$ and that it is the unique morphism that extends $\phi$. We must now check that $\Phi$ is $G$-equivariant. Consider any operation $(\theta,\delta)$ in $\SN$. We have to show by induction that, for every $g$ in $G$, 
		$$
		\Phi (g \theta, \omega_{\theta}(g)\cdot \delta)= g \Phi (\theta,\delta).
		$$
		If $\theta$ is the identity $\Id{}$, then the result is trivially verified.
		Let us suppose now by induction that $\theta=r\bigotimes_T\left(\theta_1,\ldots,\theta_{|T|}\right)$ and that the result is known for $\theta_1,\ldots, \theta_{|T|}.$
		If we follow the notation from the definition of the $G$-action on formal external norms again, then we have
		\begin{align*}
		&\textstyle \Phi(g\theta,\omega_{\theta}(g)\circ \delta) \\
        =& \textstyle
		r' \gamma \left(\bigotimes\nolimits_T^{\ma{O}},\Phi\left(h\theta_{\sigma(h)^{-1}1}\right),\ldots,\Phi\left(h\theta_{\sigma(h)^{-1}|T|}\right)\right)\cdot \left(\omega_{\theta}(g)\circ \delta\right) \\
		=& \textstyle r' \gamma\left(\bigotimes\nolimits_T^{\ma{O}},\Phi\left(h\theta_{\sigma(h)^{-1}1}\right),\ldots,\Phi\left(h\theta_{\sigma(h)^{-1}|T|}\right)\right)\cdot \left(\sigma(h)\num{\omega_1(h),\ldots, \omega_{|T|(h)}} \circ \delta\right) \\
		=& \textstyle r' \gamma \left(\left(\bigotimes\nolimits_T^{\ma{O}},\sigma(h)\right),\Phi\left(h\theta_1,\omega_{\theta_1}(h)\right),\ldots,\Phi \left(h\theta_{|T|},\omega_{\theta_{|T|}}(h)\right)\right) \cdot \delta \\
		=& \textstyle r' \gamma \left(h\bigotimes\nolimits_T^{\ma{O}},h\Phi\left(\theta_1\right),\ldots, h\Phi\left(\theta_{|T|}\right)\right)\cdot \delta \\
		=& \textstyle g\Phi \left(\theta,\delta\right). 
		\end{align*} 
	\end{proof}
\end{proposition}

We deduce from \ref{freeop} that Rubin's coherence theorem applies in our situation.
\begin{theorem}[\protect{\cite[Theorem 5.6]{Rub2}},\protect{\cite[Theorem 2.10]{Rub3}}]\label{rubco}
	The $2$-categories $\NSSlaxI$, $\NSStMI$ and $\NSSMI$ are respectively equivalent to the $2$-categories $\mathbf{AlgLax}_{\EI}$, $\mathbf{AlgStg}_{\EI}$ and $\mathbf{AlgSt}_{\EI}$ in the sense of \ref{2catalg}.
\end{theorem}

\begin{notation}
	In what follows, we denote by $\EII$ the operad $N(\EI)$ in $\gsset$.
\end{notation}

The following corollary is a direct consequence of \ref{thesubOalg}.

\begin{corollary}\label{subalgNS}
	There is a fully faithful simplicial functor 
    $$
    N\colon\NSSMI \to \Alg_{\EII} \left(\gsset\right)
    $$
    that sends an $\Ind$-normed symmetric monoidal category $\cC$ to its nerve $N(\cC)$ where the strict $2$-category $\NSSMI$ is seen as a simplicial category. 
\end{corollary}

We will now show that the operad $\EII$ is an $N_{\infty}$-operad with associated $G$-indexing system $\Ind$. We start with the following proposition, which is a reformulation in terms
of our description of $\EI$ of \cite[Theorem 4.6]{Rub1}.

\begin{proposition}\label{gromor}
	Let $H\subseteq G$ be a subgroup and $\theta$ an $H$-fixed formal external norm.
	The finite $H$-set $T_{\theta}$ (see \ref{equinorm}) is admissible . 
	\begin{proof}
		We show the result by induction on $\theta$. Note first that if $\theta$ is the identity $\Id{}$, then $T_{\Id{}}$ is the final $G$-set $*$ which is admissible by Axiom $1$ of \ref{ind}. Assume now that $\theta=r\bigotimes_T \left(\tau_1,\ldots, \tau_{|T|}\right)$ with $T$ an admissible $K$-set and $r$ a chosen coset representative of $G/K$ with $K\subseteq G$ any subgroup. It follows from the fact that $\theta$ is $H$-fixed that the coset of $G/K$ represented by $r$ is $H$-fixed and thus $H\subseteq {^rK}$. If we denote by $T'$ the $H$-set $\res{{^rK}}{H}\operatorname{c}_{r^{-1}}T$ and by $\theta'$ the formal external norm $\bigotimes_{T'}\left(\tau_1,\ldots, \tau_{|T|}\right)$, then there is an isomorphism $T_{\theta}\simeq T_{\theta'}$ of finite $H$-sets and the result holds for $\theta$ if it holds for $\theta'$. We can now assume without loss of generality that $\theta=\bigotimes_T \left(\tau_1,\ldots, \tau_{|T|}\right)$ with $T$ an admissible $H$-set. Choose a decomposition $T\simeq \bigsqcup_{i=1}^{n} H/K_i$
		and denote by $\alpha_i:\{1,\ldots, [H:K_i]\} \to \{1,\ldots,|T|\}$ the map that corresponds to the inclusion $H/K_i \hookrightarrow T$ through the ordering of $T$. If we denote by $\theta_i$ the formal external norm $\bigotimes_{H/K_i}\left(\tau_{\alpha_i(1)},\ldots,\tau_{\alpha_i\left([H:K_i]\right)}\right)$
		, then we have $T_{\theta}\simeq \bigsqcup_{i=1}^{n} T_{\theta_i}$ and using Axioms $6$ and $5$ of \ref{ind}, the result holds for $\theta$ if and only if it holds for every $\theta_j$. We can now assume without loss of generality that $T=H/K$. Since the formal external norm $\theta=\bigotimes_{H/K}\left(\tau_1,\ldots,\tau_{[H:K]}\right)$ is $H$-fixed, we have for every $h$ in $H$ that $h\tau_{\sigma(h)^{-1}j}=\tau_j$ and we deduce that $\tau_1$ is $K$-fixed and that $h_j \tau_1 = \tau_j$ with $h_j$ a representative of the $j$-th coset of $H/K$ through the chosen ordering of $H/K$. By induction, we can assume that the $K$-set $T_{\tau_1}$ is admissible and we obtain that $T_{\theta}\simeq \ind{H}{K}T_{\tau_1}$. The result now follows from Axiom $7$ of \ref{ind}.
	\end{proof}
\end{proposition}

We can now show that $\EII$ is an $N_{\infty}$-operad with associated $G$-indexing system $\Ind$.

\begin{proposition}
	The operad $\EII$ in $\gsset$ is an $N_{\infty}$-operad with $\Ind$ as its $G$-indexing system.
	\begin{proof}
		Let $\Gamma \subseteq G\times \Sigma_n$ be any subgroup.
		The functors $\widetilde{(-)}$ and $N$ preserve limits and it follows that
		$$
		{(\EII)_n}^{\Gamma}\simeq N\left (\widetilde{{\SN_n}^{\Gamma}}\right).
		$$
		For every $G$-set $X$, the simplicial set $N(\widetilde{X})$ is respectively contractible or empty if and only if $X$ is non-empty or empty and we deduce that ${(\EII)_n}^{\Gamma}$ is respectively contractible or empty if and only if ${\SN_n}^{\Gamma}$ is non-empty or empty.
		The freeness of the symmetric actions on $\EII$ now follows from the freeness of the symmetric actions on $\SN$ by considering subgroups of $G\times \Sigma_n$ of the form $\{1_G\}\times \Lambda$ with $\Lambda\subseteq \Sigma_n$.
		To conclude that $\EII$ is an $N_{\infty}$-operad with $\Ind$ as $G$-indexing system, it is now enough to show that, for every finite $H$-set $T$, ${\SN_n}^{\Gamma_T}$ is not empty if and only if $T$ is admissible. Assume first that $T$ is admissible. The $n$-ary operation $\bigotimes_T$ of $\SN$ is clearly fixed by $\Gamma_T$ and ${\SN_n}^{\Gamma_T}$ is not empty. Suppose now that $(\theta,\delta)$ is a $n$-ary operation of $\SN$ fixed by $\Gamma_T$. If we denote by $\sigma:H \to \Sigma_{n}$ the homomorphism corresponding to $T$ through the chosen ordering $T\simeq \{1,\ldots,|T|\}$, we deduce that $\theta$ is $H$-fixed and that $\omega_{\theta}(h)\circ \delta =\delta \circ \sigma(h)$. The permutation $\delta$ can now be seen as an $H$-equivariant isomorphism between $T_{\theta}$ and $T$ and, by \ref{gromor} and Axiom $2$ of \ref{ind}, $T$ is admissible. The result follows.
	\end{proof}
\end{proposition}

\section{Mackey functors vs \texorpdfstring{$\EII$-algebras}{EI-algebras}}\label{guimay}
In this section, we define the $\Ind$-effective Burnside category $\SpaG$ and we use the $2$-categorical unfurling of \cite[Proposition $2.24$]{TOB1} to construct a fully faithful functor $\fG :\SpaG \to N_{\Delta}\left(\NSSMI_{(2,1)}\right)$ with $N_{\Delta}$ the homotopy coherent nerve and $\NSSMI_{(2,1)}$ the underlying $(2,1)$-category of $\NSSMI$. We conclude the section by showing that $\cG$ induces a functor $\theta:\SpaG \to \Ninf$ whose left Kan extension $\Theta:\func_{\times}(\opo{\SpaG},\sS) \to \Ninf$ along the Yoneda embedding $\SpaG \to \func_{\times}(\opo{\SpaG},\sS)$ is an equivalence of $\infty$-categories. This will prove \ref{THETH} in \ref{subsecmainth}.

\subsection{Orbital subcategories and the \texorpdfstring{$\Ind$}{I}-effective Burnside category}
In this subsection, we explain how the $G$-indexing system $\Ind$ leads to the definition of the $\Ind$-effective Burnside category $\SpaG$. 
\begin{definition}[\protect{\cite[Definition $2.5$]{dede}}]\label{finitecopchap2}
	Let $\mathscr{O}$ be a small $\infty$-category. The \emph{finite cocompletion of $\mathscr{O}$} denoted by $\FFF{\mathscr{O}}$ is the full subcategory of $\Psh (\mathscr{O})$ spanned by finite coproducts of representables. The $\infty$-category $\FFF{\mathscr{O}}$ admits finite coproducts and satisfies the following universal property: for every $\infty$-category $\cC$ that admits finite coproducts, the natural functor $\mathscr{O} \to \FFF{\mathscr{O}}$ induces an equivalence of $\infty$-categories
	$$
	\fun_{\sqcup}(\FFF{\mathscr{O}},\cC)\to \func (\mathscr{O},\cC)
	$$
	where $\func_{\sqcup}(\FFF{\mathscr{O}},\cC)$ is the full subcategory of $\func (\FFF{\mathscr{O}},\cC)$ of those functors $\phi:\FFF{\mathscr{O}}\to \cC$ that preserve finite coproducts.
\end{definition}

\begin{example}
	The finite cocompletion of $\OG$ is equivalent to the category $\gset$.
\end{example}

We recall now from \cite{CLL1} the notion of orbital subcategories.

\begin{definition}[\protect{\cite[Definition $4.2.2$]{CLL1}}]\label{orb}
	Let $\mathscr{O}$ be a small $\infty$-category. A wide subcategory $\mathscr{T}\subseteq \mathscr{O}$ is \emph{orbital} if the base change of a morphism in $\FFF{\mathscr{T}}$ along any morphism in $\FFF{\mathscr{O}}$ exists and belongs to $\FFF{\mathscr{T}}$. Equivalently, for every pullback diagram 
	$$
	\begin{tikzcd}
	A \arrow[r,"u"] \arrow[d,"v"'] & B \arrow[d,"w"] \\
	C \arrow[r,"z"'] & D
	\end{tikzcd}
	$$
	in $\Psh (\mathscr{O})$ with $B$, $C$ and $D$ in $\mathscr{O}$, we can decompose $A$ as a disjoint union $\bigsqcup_{i=1}^n A_i$ of objects of $\mathscr{O}$ and, if we write $v$ through such decomposition as $(v_i)_{i=1}^n:\bigsqcup_{i=1}^n A_i \to C$, then $v_i$ belongs to $\mathscr{T}$. A small $\infty$-category $\cC$ is \emph{orbital} if it is orbital when regarded as a subcategory of itself.
\end{definition}

\begin{example}\label{orbJ}
	For every $G$-indexing system $\ma{J}$, consider the subcategory $\mathscr{O}_{\ma{J}}$ of the orbit category $\OG$ spanned by the morphisms $u:G/K\to G/H$ such that $H/K^g$ belongs to $\ma{J} (H)$, where $g$ is a representative of the $K$-fixed coset of $G/H$ that corresponds to $u$. The category $\sO_{\cJ}$ is an orbital subcategory of $\sO_G$.
\end{example}

\begin{definition}\label{bubu}
	If $\mathscr{O}$ is an $\infty$-category and $\mathscr{T}$ an orbital subcategory of $\mathscr{O}$, then we denote by $\spa{\mathscr{O}}{\mathscr{T}}$ the effective Burnside category $A^{\emph{eff}}\left(\FFF{\mathscr{O}},\FFF{\mathscr{O}},\FFF{\mathscr{T}}\right)$ in the sense of \cite[Definition 5.10]{Bar1}.
\end{definition}

The following proposition is a direct consequence of \cite[Lemma C.3]{BaHo1}.
\begin{proposition}\label{addispa}
	Let $\mathscr{O}$ be an $\infty$-category and $\mathscr{T}$ an orbital subcategory of $\mathscr{O}$. The $\infty$-category $\spa{\mathscr{O}}{\mathscr{T}}$ is semi-additive.
\end{proposition}

We now prove that the poset of $G$-indexing systems, and equivalently the poset of $N_{\infty}$-operads, is equivalent to the poset of orbital subcategories of the orbit category $\OG$.

\begin{proposition}\label{orbind}
	Denote by $\OrbG$ the poset of orbital subcategories of $\OG$ under inclusion. The morphism of posets $\IndG \to \OrbG$ that sends $\ma{J}$ to the orbital subcategory $\mathscr{O}_{\ma{J}}$ of \ref{orbJ} is an isomorphism.
	
	\begin{proof}
		The category $\gset$ is equivalent to the finite cocompletion of $\OG$ and there is an isomorphism of posets between $\OrbG$ and the poset of wide, pullback stable, finite coproduct complete subcategories of $\gset$ that sends an orbital subcategory $\mathscr{T}\subseteq \OG$ to its finite cocompletion $\FFF{\mathscr{T}}$ seen as a wide subcategory of $\gset$. The result now follows from \cite[Theorem 1.4]{MR3773736}.
	\end{proof}
\end{proposition}

\begin{definition}
	The \emph{$\Ind$-effective Burnside category} $\SpaG$ is the $\infty$-category $\spa{\OG}{\mathscr{O}_{\ma{I}}}$ in the sense of \ref{bubu}. A $G$-equivariant map $u:A\to B$ between finite $G$-sets is \emph{admissible} if it belongs to $\mathscr{O}_{\ma{I}}$. 
\end{definition}

The following proposition gives a simple way to characterise admissible maps.

\begin{proposition}\label{adm}
	Let $u:A\to B$ be a $G$-equivariant map between finite $G$-sets. If we choose a decomposition $B\simeq \bigsqcup_{i=1}^n G/H_i$, then the map $u$ is admissible if and only if for every $i=1,\ldots,n$, the finite $H_i$-set $u^{-1}(1_G H_i)$ is admissible.
	\begin{proof}
		Choose a decomposition $B\simeq \bigsqcup_{i=1}^n G/H_i$ and write $u$ as $\bigsqcup_{i=1}^n u_i$ with $u_i:A_i \to G/H_i$ where $A_i=u^{-1}(G/H_i)$. By definition, the result stands for $u$ if and only if it holds for every $u_i$ and we can assume without loss of generality that $B=G/H$.
		Choose a decomposition $A\simeq \bigsqcup_{j=1}^m G/K_j$ and write $u$ through this decomposition as $(u_j)_{j=1}^m$ with $u_j:G/K_j \to G/H$. By definition, $u$ is admissible if and only if every $u_j$ is admissible. Moreover, $u^{-1}(1_G H)\simeq \bigsqcup_{j=1}^m u_j^{-1}(1_G H)$, and if follows from Axioms $5$ and $6$ of \ref{ind} that the result stands for $u$ if and only if it holds for every $u_i$ so we can assume without loss of generality that $A=G/K$. If $u:G/K\to G/H$ is now given by a $K$-fixed coset representative $g$ of $G/H$, the result follows from the fact that $u^{-1}(1_GH)$ is exactly $H/K^g$, which is by definition an admissible $H$-set if and only if $u$ is admissible.
	\end{proof}
\end{proposition}

The $\Ind$-effective Burnside category $\SpaG$ induces a notion of $\Ind$-Mackey functors.

\begin{definition}
	An \emph{$\Ind$-Mackey functor} is a functor $\ma{M}:\opo{\SpaG}\to \sS$ that preserves finite products. We denote by $\Mack$ the $\infty$-category of $\Ind$-Mackey functors $\func_{\times}(\opo{\SpaG},\sS)$.
\end{definition}

\begin{remark}
	By \ref{addispa}, the $\infty$-category $\SpaG$ is semi-additive and it follows from \cite[Corollary 2.5]{MR3450758} that the $\infty$-category $\Mack$ is equivalent to the $\infty$-category $\func_{\oplus}(\SpaG,\operatorname{Mon}_{E_{\infty}}(\sS))$ with $\operatorname{Mon}_{E_{\infty}}(\sS)$ the $\infty$-category of $E_{\infty}$-monoids in $\sS$. The $\infty$-category $\Mack$ is thus a natural $\infty$-categorical generalisation of the category of classical $\Ind$-Mackey functors with values in commutative monoids. 
\end{remark}

\subsection{Unfurling for 2-categories and orbital subcategories}
In this subsection, we recall the $2$-categorical unfurling of \cite[Proposition 2.24]{TOB1} in the case of $\spa{\mathscr{O}}{\mathscr{T}}$ with $\mathscr{O}$ a category and $\mathscr{T}\subseteq \mathscr{O}$ an orbital subcategory of $\mathscr{O}$. 

\begin{proposition}[\protect{\cite[Proposition 2.22]{TOB1}}]\label{univprop}
	Let $\cC$ be a strict $2$-category and let $\phi:\FFF{\mathscr{O}} \to \cC$ be a strict $2$-functor such that for every morphism $u:A\to B$ in $\FFF{\mathscr{T}}$, the $1$-cell $u_!:=\phi(u)$ admits a right adjoint $u^*$ satisfying the Beck-Chevalley condition: for every pullback square in $\FFF{\mathscr{O}}$ as on the left in
	$$
	\begin{tikzcd}
	A\arrow[r,"u"] \arrow[d,"v"']  & B \arrow[d,"w"] &\text{} & \phi(A)  \arrow[dr,Rightarrow,shorten >=4mm,shorten <=4mm] \arrow[r,"u_!"] & \phi(B)  \\
	C \arrow[r,"z"'] &D &\text{} & \phi(C) \arrow[r,"z_!"'] \arrow[u,"v^*"] & \phi(D) \arrow[u,"w^*"'] 
	\end{tikzcd}
	$$
	with vertical arrows in $\FFF{\mathscr{T}}$, the canonical mate depicted on the right is an isomorphism.
	Then there is a unique functor $\Phi:\spa{\mathscr{O}}{\mathscr{T}}\to N_{\Delta}\left(\cC_{(2,1)}\right)$ that sends a $2$-simplex of the form 
	$$
	\begin{tikzcd}[column sep=small]
	&& A \arrow[rd,"j"] \arrow[ld,"i"']&& \\
	&B \arrow[rd,"l"'] \arrow[ld,"k"'] &&C \arrow[rd,"n"] \arrow[ld,"m"]& \\
	D&&E&&F&
	\end{tikzcd}
	$$
	to the $2$-simplex 
	$$
	\begin{tikzcd}[column sep=small]
	&|[alias=Z]| \phi (E)   \arrow[rd,"n_!m^*"]  & \\
	\phi(D) \arrow[ru,"l_! k^*"] \arrow[rr,"(nj)_!(ki)^*"'name=O] &  &\phi(F) \arrow[from=O,to=Z,Rightarrow,shorten >=2mm,shorten <=3mm]
	\end{tikzcd}
	$$
	given by the pasting 
	$$
	\begin{tikzcd}[column sep=small]
	&  & A \arrow[dd,Rightarrow,shorten >=7mm,shorten <=7mm]  \arrow[rd,"j_!"] && \\
	&|[alias=B]| B \arrow[rd,"l_!"'] \arrow[ru,"i^*"]  & &C \arrow[rd,"n_!"]  & \\
	\arrow[ru,"k^*"] \arrow[rruu, bend left=1.8cm,"(ki)^*"name=A] D&&E\arrow[ru,"m^*"']&&F.& \arrow[from=A,to=B,Rightarrow, shorten <=2.3mm]
	\end{tikzcd}
	$$
	Here $\cC_{(2,1)}$ is the $(2,1)$-category obtained from $\cC$ by throwing away the non-invertible $2$-cells.    
\end{proposition}

\begin{definition}\label{Ichev}
	Let $\cC$ be a strict $2$-category. A strict $2$-functor $\phi:\FFF{\mathscr{O}} \to \cC$ satisfies the \emph{$\mathscr{T}$-Beck Chevalley condition} if it satisfies the conditions of \ref{univprop}. In the case where $\mathscr{O}$ is $\OG$ and where $\mathscr{T}$ is the orbital subcategory $\mathscr{O}_{\Ind}$ of $\OG$ that corresponds to $\Ind$, we say that $\phi$ satisfies the \emph{$\Ind$-Beck-Chevalley condition}.
\end{definition}

\begin{notation}
	Let $\cD$ be a category with finite coproducts. If $A_1, \ldots , A_n$ and $B_1,\ldots ,B_m$ are objects of $\cD$, for every $i=1,\ldots, n$, $f_i:A \to B_{\alpha(i)}$ is a morphism in $\cD$ and $\alpha:\{1,\ldots n\}\to \{1,\ldots,m\}$ is a map, then we denote by $\left((f_i)_{i=1}^n,\alpha\right):\bigsqcup_{i=1}^n A_i \to \bigsqcup_{j=1}^m B_j$ the morphism $(\iota_{\alpha(i)} \circ f_i)_{i=1}^n$ with $\iota_j:B_j \to \bigsqcup_{j=1}^m B_j$ the canonical inclusion. Note that if $\mathscr{E}$ is a small category, then every morphism in $\FFF{\mathscr{E}}$ can be written canonically as $((f_i)_{i=1}^n,\alpha):\bigsqcup_{i=1}^n A_i \to \bigsqcup_{j=1}^m B_j$ with $A_i$ and $B_j$ in $\mathscr{E}$.
\end{notation}

In the case where $\cC$ is a semi-additive strict $2$-category, we now give a simple criterion to check whether a strict $2$-functor $\phi:\FFF{\mathscr{O}}\to \cC$ that preserves finite coproducts satisfies the $\mathscr{T}$-Beck-Chevalley condition.

\begin{proposition}\label{univquasi}
	Let $\cC$ be a semi-additive strict $2$-category and $\phi:\FFF{\mathscr{O}}\to \cC$ a strict $2$-functor that preserves finite coproducts. The strict $2$-functor $\phi:\FFF{\mathscr{O}}\to \cC$ satisfies the $\mathscr{T}$-Beck-Chevalley condition if and only if for every morphism $u:A\to B$ in $\FFF{\mathscr{T}}$ with $B$ in $\mathscr{O}$, the $1$-cell $\phi(u):= u_!$ admits a right adjoint $u^*$ such that for every pullback square in $\FFF{\mathscr{O}}$ as on the left in 
	$$
	\begin{tikzcd}
	A \arrow[d,"v"'] \arrow[r,"u"]  & B \arrow[d,"w"] &\text{} &\arrow[dr,Rightarrow,shorten >=4mm,shorten <=4mm]  \phi (A) \arrow[r,"u_!"] & \phi(B)  \\
	C \arrow[r,"z"'] &D &\text{}  & \phi(C) \arrow[u,"v^*"] \arrow[r,"z_!"'] & \phi(D) \arrow[u,"w^*"']
	\end{tikzcd}
	$$
	with $w$ in $\FFF{\mathscr{T}}$, $z$ in $\mathscr{O}$, and $C$ and $D$ in $\mathscr{O}$, the canonical mate depicted on the right is an isomorphism.  Moreover, if a strict $2$-functor $\phi:\FFF{\mathscr{O}}\to \cC$ that preserves finite coproducts satisfies the $\mathscr{T}$-Beck-Chevalley condition, the induced functor $\Phi:\spa{\mathscr{O}}{\mathscr{T}}\to N_{\Delta}\left(\cC_{(2,1)}\right)$ is semi-additive.
	\begin{proof}
		It is clear that the $\mathscr{T}$-Beck-Chevalley condition implies the second condition of the proposition. Assume now that $\phi$ satisfies the second condition. Let $u:A \to B$ be a morphism in $\FFF{\mathscr{T}}$. Choose a decomposition $B\simeq \bigsqcup_{k=1}^p B_k$ with $B_k$ in $\mathscr{O}$ through which $u$ can be written as $\bigsqcup_{k=1}^p u_k$ with $u_k:A_k \to B_k$ in $\FFF{\mathscr{T}}$. By assumptions, $(u_k)_!$ admits a right adjoint $(u_k)^*$ and it follows that $u_!$ also admits a right adjoint given by $u^*=\bigoplus_{k=1}^p (u_k)^*$. Consider now any pullback square
		\begin{equation}\label{squ}\tag{S}
		\begin{tikzcd}
		A \arrow[r,"u"] \arrow[d,"v"'] & B \arrow[d,"w"] \\
		C \arrow[r,"z"'] & D
		\end{tikzcd}
		\end{equation}
		in $\FFF{\mathscr{O}}$ with $v$ and $w$ in $\FFF{\mathscr{T}}$. Choose decompositions $C \simeq \bigsqcup_{i=1}^n C_i$ and $D\simeq \bigsqcup_{j=1}^m D_j$ with $C_i$ and $D_j$ in $\mathscr{O}$ through which we have  $z=((z_i)_{i=1}^n,\alpha)$ and $w=\bigsqcup_{j=1}^m w_j$ with $w_j:B_j \to D_j$. For every $i=1,\ldots, n$ , consider the following pullback square 
		\begin{equation}\label{squa}\tag{$S_{i}$}
		\begin{tikzcd}
		A_{i} \arrow[r,"u_{i}"] \arrow[d,"v_{i}"'] & B_{\alpha(i)} \arrow[d,"w_{\alpha(i)}"] \\
		C_i \arrow[r,"z_i"'] & D_{\alpha(i)}.   
		\end{tikzcd}
		\end{equation}
		By disjointness and universality of finite sums in $\FFF{\mathscr{O}}$, we have an isomorphism $A\simeq \bigsqcup_{i=1}^n A_{i}$ through which we can write $u$ and $v$ respectively as $((u_i)_{i=1}^n,\alpha):\bigsqcup_{i=1}^n A_i \to  \bigsqcup_{j=1}^m B_j$ and $\bigsqcup_{i=1}^n v_{i}$. The morphisms $v_{i}$ and $w_{j}$ both belong in $\FFF{\mathscr{T}}$ and $(v_{i})_!$ and $(w_j)_!$ respectively admit right adjoints $(v_{i})^*$ and $(w_j)^*$ by assumption. Moreover, the mate $\omega_{i}$ depicted in the following diagram associated with the square \ref{squa} is an isomorphism by assumption
		$$
		\begin{tikzcd}
		\phi \left(A_{i}\right) \arrow[dr,Rightarrow,shorten >=4.5mm,shorten <=4.5mm,"\omega_{i}"] \arrow[r,"(u_{i})_!"]  & \phi \left (B_{\alpha(i)}\right)  \\
		\phi(C_i) \arrow[r,"(z_i)_!"'] \arrow[u,"(v_{i})^*"] & \arrow[u,"(w_{\alpha(i)})^*"']\phi \left(D_{\alpha(i)}\right).   
		\end{tikzcd}
		$$
		The morphisms $w_!=\bigoplus_{j=1}^m (w_j)_!$ and $v_!=\bigoplus_{i=1}^n (v_i)_!$ respectively admit right adjoints $w^*=\bigoplus_{j=1}^m (w_j)^*$ and $v^*=\bigoplus_{i=1}^n (v_i)^*$ and
		by direct inspection, the mate $\omega$ depicted in the following diagram associated with the square \ref{squ}
		$$
		\begin{tikzcd}
		\varphi(A) \arrow[dr,Rightarrow,shorten >=4mm,shorten <=4mm,"\omega"] \arrow[r,"u_!"]  & \phi(B)  \\
		\phi(C) \arrow[r,"z_!"'] \arrow[u,"v^*"] & \arrow[u,"w^*"']\phi(D)   
		\end{tikzcd}
		$$
		corresponds to $(\omega_{i,j})_{i,j}$ through the equivalence 
		$$
		\cC (\phi(C),\phi(B)) \simeq \prod_{i=1}^n\prod_{j=1}^m \cC(\phi(C_i),\phi(B_j)),
		$$
		where we define $\omega_{i,j}$ by $\Id{0_{\phi(C_j),\phi(B_i)}}$ whenever $\alpha(i)\neq j$. Every $\omega_{i,j}$ is an isomorphism and it follows that $\omega$ is an isomorphism and that $\phi$ satisfies the $\mathscr{T}$-Beck-Chevalley condition. If $\phi:\FFF{\mathscr{O}}\to \cC$ is a strict $2$-functor that preserves finite coproducts, then it follows directly from \cite[Proposition 4.3]{Bar1} that the induced functor $\Phi :\spa{\mathscr{O}}{\mathscr{T}}\to N_{\Delta} \left(\cC_{(2,1)}\right)$ is semi-additive.
	\end{proof}
\end{proposition}

\subsection{From \texorpdfstring{$\SpaG$}{G-spaces} to \texorpdfstring{$\NSSMI$}{NSSMI}}\label{exa}

In this subsection, we construct a strict $2$-functor $\fG:\gset \to \NSSMI$ that satisfies the $\Ind$-Beck Chevalley condition in the sense of \ref{Ichev}.

\begin{notation}
	If $n\ge 0$, we denote by $\num{n}$ the finite set $\{1,\ldots, n\}$.
\end{notation}

\begin{definition}\label{theGcat}
	Let $A$ be a finite $G$-set. The $G$-category $\fG(A)$ has as objects pairs $(\theta, \cL)$ with $\theta$ a formal external norm and $\cL=\left(a_1,\ldots, a_{|\theta|}\right)$ a sequence of elements of $A$. A morphism $\alpha:\left(\theta,\left(a_1,\ldots,a_{|\theta|}\right)\right) \to \left(\tau,\left(b_1,\ldots,b_{|\tau|}\right)\right)$ in $\fG(A)$ is a map $\alpha:\num{|\theta|}\to \num{|\tau|}$ such that $a_i=b_{\alpha(i)}$ for every $i$ in $\num{|\theta|}$. The action of $G$ on $\fG(A)$ is defined on objects by 
	$$
	g\cdot \left(\theta,\left(a_1,\ldots, a_{|\theta|}\right)\right)=\left(g\theta,\left(g a_{{\omega_{\theta}(g)}^{-1}1},\ldots ,g a_{{\omega_{\theta}(g)^{-1}|\theta|}}\right)\right)
	$$
	and on a morphism $\alpha:\left(\theta,\left(a_1,\ldots,a_{|\theta|}\right)\right) \to \left(\tau,\left(b_1,\ldots,b_{|\tau|}\right)\right)$ by $g\cdot \alpha=\omega_{\tau}(g)\circ \alpha \circ \omega_{\theta}(g)^{-1}$. It follows from \ref{actionwelldef} that these $G$-actions are well-defined. The $G$-category $\fG(A)$ is an $\Ind$-normed symmetric monoidal category with the following structure:
	\begin{enumerate}
		\item the unit $e$ is the object $\left(\bigotimes_{\emptyset},()\right)$
		with $()$ the empty list;
		\item the tensor product $\otimes:\fG(A)\times \fG(A) \to \fG(A)$ is defined on objects by
		$$
		\left(\theta,\cL_1\right) \otimes (\tau,\cL_2)=\left(\theta \otimes \tau,\cL_1\sqcup \cL_2\right)
		$$
		and on morphisms by $\alpha \otimes \beta =\alpha \sqcup \beta$ where $\theta \otimes \tau$ is the formal external norm given by $\gamma \left(\bigotimes_{\num{2}},\theta,\tau\right)$ with $\num{2}$ seen as a trivial $G$-set;
		\item given an admissible $H$-set $T$, we define the external $T$-norm $\bigotimes_T :\fG(A)^{\times T}\to \fG(A)$ on objects by
		$$
		\textstyle \bigotimes\nolimits_T \left((\theta_1,\cL_1),\ldots, \left(\theta_{|T|},\cL_{|T|}\right)\right)=
		\left(\bigotimes\nolimits_T\left(\theta_1,\ldots, \theta_{|T|}\right),\bigsqcup_{i=1}^{|T|}\cL_i\right).
		$$
		and on morphisms by $\bigotimes_T \left(\alpha_1,\ldots,\alpha_{|T|}\right)=\bigsqcup_{i=1}^{|T|}\alpha_i$.
	\end{enumerate}
	If $X=\left(\theta,\left(a_1,\ldots, a_{|\theta|}\right)\right)$ is an object of $\fG(A)$, then the \emph{length} of $X$ is the length of $\theta$ and we denote it by $|X|$.
	If $X$, $Y$ and $Z$ are objects of $\fG(A)$ with respective lengths $n$, $m$ and $k$, the natural equivariant isomorphisms $\alpha_{X,Y,Z}:(X \otimes Y) \otimes Z \to X \otimes (Y\otimes Z)$, $\lambda_X:e\otimes X \to X$, $\rho_X:X \otimes e\to X$ and $\beta_{X,Y}:X \otimes Y \to Y \otimes X$ are respectively given by the permutations $1_{\Sigma_{n+m+k}}$, $1_{\Sigma_n}$, $1_{\Sigma_n}$ and $\tau (n,m)$ with $\tau (n,m)$ the block permutation associated with the unique non trivial permutation $\tau$ of $\Sigma_2$. If $X_1, \ldots, X_{|T|}$ are objects of $\fG(A)$ with respective lengths $n_1,\ldots, n_{|T|}$, the natural nonequivariant isomorphism $v_{T}:\bigotimes_{T}\left(X_1,\ldots,X_{|T|}\right)\to \bigotimes_{|T|}\left(X_1,\ldots,X_{|T|}\right)$ is given by the permutation $1_{\Sigma_{n_1+\cdots +n_{|T|}}}$. 
\end{definition}

\begin{remark}
	If $A$ is a finite $G$-set, an object $\left(\theta,\left(a_1,\ldots,a_{|\theta|}\right)\right)$ of $\fG(A)$ can be seen as the planar rooted tree associated with $\theta$ where the $i$-th leaf represents $a_i$.
	\begin{center}
		\begin{tikzpicture}[grow via three points={%
			one child at (0,1) and two children at (-0.75,1) and (0.75,1)}]
		\node at (0,-4.5) {$r_1\bigotimes_{T_1}$} child{node {$r_2\bigotimes_{T_2}$}child{node{$r_4\bigotimes_{ \emptyset}$} } child{node {$a_1$}} child{node {$a_2$}}} child {node {$\cdots$}} child{node {$\cdots$}} child{node {$r_3\bigotimes_{T_3}$} child{ node {$r_5\bigotimes_{\emptyset}$}} child{ node {$a_{|\theta|}$}}} ;
		\end{tikzpicture}
	\end{center}
\end{remark}

\begin{notation}
	In what follows, given an element $a$ of $A$, we will denote by $a$ the object $(\Id{},(a))$ of $\fG(A)$. This gives a canonical inclusion $A\hookrightarrow \fG(A)$. 
\end{notation}

The construction of \ref{theGcat} is functorial. For every $G$-equivariant map $u:A\to B$, the definition of the strict $\Ind$-normed functor $\fG (u):\fG (A)\to \fG(B)$ is part of the following construction.

\begin{construction}\label{thecons}
	Let $A$ and $B$ be two finite $G$-sets and $\phi:A\to \fG(B)$ a $G$-functor with $A$ seen as a discrete $G$-category. We define a strict $\Ind$-normed functor $\conss(\phi):\fG(A)\to \fG(B)$ as follows: 
	on an object $X=\left(\theta,\left(a_1,\ldots,a_{|\theta|}\right)\right)$ of $\fG(A)$, the strict $\Ind$-normed functor $\conss(\phi)$ is defined by the following recursion: 
	\begin{enumerate}
		\item if $\theta=\Id{}$, then $\conss(\phi)(X)=\phi(a)$;
		\item if $\theta=r\bigotimes_T \left(\theta_1,\ldots, \theta_{|T|}\right)$, then $\conss(\phi) (X)$ is
		$$\textstyle
		r\bigotimes\nolimits_T\left(\conss(\phi)\left(\theta_1,\left(a_1,\ldots,a_{|\theta_1|}\right)\right),\ldots, \conss(\phi)\left(\theta_{|T|},\left(a_{|\theta|-|\theta_{|T|}|},\ldots,a_{|\theta|}\right)\right)\right)
		$$
	\end{enumerate} 
	and on a morphism $\alpha:\left(\theta,\left(a_1,\ldots,a_{|\theta|}\right)\right)\to \left(\tau,\left(b_1,\ldots,b_{|\tau|}\right)\right)$, the map $\conss(\phi) (\alpha)$ is defined by the block map 
	$$
	\alpha\left(k_1,\ldots, k_{|\theta|}\right)\left(l_1,\ldots, l_{|\tau|}\right):\num{k_1+\cdots +k_{|\theta|}}\to \num{l_1+\cdots+l_{|\tau|}}
	$$
	that sends the copy of $\num{k_i}$ included in $\num{k_1}\sqcup \cdots \sqcup \num{k_{\theta}} \simeq \num{k_1+\cdots +k_{|\theta|}}$ to the copy of $\num{k_i} =\num{l_{\alpha(i)}}$ included in  $\num{l_1}\sqcup \cdots \sqcup \num{l_{|\tau|}}  \simeq \num{l_1+\cdots +l_{|\tau|}}$
	where $k_i$ and $l_j$ are respectively the length of the image of $a_i$ and $b_j$ by $\phi$ (we use here that $a_i=b_{\alpha(i)}$). This construction is functorial: if $\phi,\psi:A \rightrightarrows \fG(B)$ are two $G$-functors and $\omega:\phi \Rightarrow \psi$ a $G$-natural transformation  between them, then we obtain an $\Ind$-normed natural transformation $\conss (\omega) :\conss (\phi ) \Rightarrow \conss (\psi)$ defined by
	\[
	\conss (\omega)_{\left(\theta,\left(a_1,\ldots,a_{|\theta|}\right)\right)}=\bigsqcup_{i=1}^{|\theta|}\omega_{a_i}.
	\]
	We denote by $\conss:\func_G (A,\fG(B)) \to \NSSMI (\fG(A),\fG(B))$ the functor provided by this construction.
\end{construction}

\begin{lemma}
	For finite $G$-sets $A$ and $B$, the functor 
	$$
	\conss:\func_G (A,\fG(B)) \to \NSSMI (\fG(A),\fG(B))
	$$
	is fully faithful.
	\begin{proof}
		By definition, if $\phi,\psi :A \rightrightarrows \fG(B)$ are two $G$-functors, then, because $\conss (\phi)$ and $\conss (\psi)$ are strict $\Ind$-normed functors, an $\Ind$-normed natural transformation $\Omega:\conss (\phi) \Rightarrow \conss (\psi)$ is determined by its values on elements of $A$ and $\Omega$ must be equal to $\conss(\omega)$ where $\omega$ is the natural transformation $\omega :\phi \Rightarrow \psi$ obtained by restricting $\Omega$ to $A$. The result follows.
	\end{proof}
\end{lemma}

\begin{proposition}\label{composi}
	The following diagram of functors is strictly commutative.
	$$
	\begin{tikzcd}[font=\footnotesize]
	\func_G (A,\fG(B)) \times \func_G(B,\fG(C)) \arrow[r,"\Id{\func_G (A,\fG(B))}\times \conss"] \arrow[dd,"\conss \times \conss"'] & \func_G(A,\fG(B)) \times \NSSMI (\fG(B),\fG(C)) \arrow[d,"\circ"] \\
	\text{} & \func_G(A,\fG(C)) \arrow[d,"\conss"] \\
	\NSSMI (\fG(A),\fG(B))\times \NSSMI (\fG(B),\fG(C)) \arrow[r,"\circ"'] & \NSSMI (\fG(A),\fG(B))  
	\end{tikzcd}
	$$
	\begin{proof}
		Let $\phi:A\to \fG(B)$ and $\psi:B\to \fG(C)$ be two $G$-functors. It follows from a direct induction on formal external norms
		that the strict $\Ind$-normed functors $\conss (\psi)\circ \conss (\phi)$ and $\conss (\conss(\psi)\circ \phi )$ agree on objects. Consider now a morphism $\alpha:X \to Y$ in $\fG(A)$ between two objects $X=\left(\theta,\left(a_1,\ldots, a_{|\theta|}\right)\right)$ and $Y=\left(\tau,\left(b_1,\ldots,b_{|\tau|}\right)\right)$ of $\fG(A)$. If we denote by $X_{i}=\left(\lambda_i,\left(c_1^i,\ldots, c^i_{|\lambda_i|}\right)\right)$ the image of $a_i$ by $\phi$, then by definition, the functor $\conss (\psi)\circ \conss (\phi)$ sends $\alpha$ on the block map $\alpha \left(|\lambda_1|,\ldots,|\lambda_{\theta}|\right)\left(k_1^1,\ldots,k_{|\lambda_1|}^1,\ldots, k_{\lambda_{|\theta|}}^{|\theta|}\right)$ where $k_{j}^i$ is the length of $\psi (c^i_j)$. It follows that the image of $\alpha$ under $\conss (\psi)\circ \conss (\phi)$ is equal to the block map $\alpha\left(k_1^1+ \cdots +k_{|\lambda_1|}^1, \ldots , k_{1}^{|\theta|}+ \cdots + k_{\lambda_{|\theta|}}^{|\theta|}\right)$ which is also the image of $\alpha$ under $\conss (\conss(\psi)\circ \phi )$. A similar argument shows that the diagram of the proposition commutes on morphisms.
	\end{proof}
\end{proposition}

If $u:A\to B$ is a $G$-equivariant map, then we denote by $\fG(u):\fG(A)\to \fG(B)$ the image under $\conss$ of the composition $A\overset{u}{\to} B \to \fG(B)$ and it follows from \ref{composi} that we get a strict $2$-functor $\fG:\gset \to \NSSMI$.
The $2$-category $\NSStMI$ is semi-additive by \ref{2sa} and we can thus prove the following proposition.
\begin{proposition}\label{papro}
	The functor $\fG$, seen as a  strict $2$-functor $\fG:\gset \to \NSStMI$ through the inclusion $\NSSMI \hookrightarrow \NSStMI$, preserves finite coproducts.
	\begin{proof}
		By \ref{2sa}, the $2$-category $\NSStMI$ is semi-additive and, given two finite $G$-sets $A$ and $B$, the coproduct of $\fG(A)$ and $\fG(B)$ in $\NSStMI$ is given by their product $\fG(A)\times \fG(B)$.
		We have to show now that the canonical strong $\Ind$-normed functor 
		$$
		\Phi=\left(\fG(i_A) \circ p_A\right) \otimes (\fG(i_B)\circ p_B):\fG(A)\times \fG(B) \to \fG(A\sqcup B)
		$$ 
		is an equivalence, where $i_A:A\to A\sqcup B$, $i_B:B\to A\sqcup B$, $p_A:\fG(A)\times \fG(B) \to \fG(A)$ and $p_B:\fG(A)\times \fG(B) \to \fG(B)$ are respectively the canonical inclusions and projections.
		The functor $\Phi$ sends an object $\left(\left(\theta,\left(a_1,\ldots,a_{|\theta|}\right)\right),\left(\tau,\left(b_1,\ldots, b_{|\tau|}\right)\right)\right)$ of $\fG(A)\times \fG(B)$ to $\left( \theta \otimes \tau, \left(a_1,\ldots, a_{\theta},b_1,\ldots, b_{|\tau|}\right)\right)$ and a morphism $\left(\alpha,\beta\right)$ of $\fG(A)\times \fG(B)$ to $\alpha \otimes \beta$. It admits an inverse in $\NSStMI$ that can be constructed as follows: consider the $G$-functor $\psi_A: A \sqcup B\to \fG(A)$ given respectively by the canonical inclusion $A \to \fG(A)$ and by the $G$-functor $B\to \fG(A)$ constant on the unit $e$. Denote by $\Psi_A:\fG(A\sqcup B)\to \fG(A)$ the strict $\Ind$-normed functor $\conss (\psi_A)$. Symmetrically, we obtain a strict $\Ind$-normed functor $\Psi_B:\fG (A\sqcup B)\to \fG(B)$ and we get a strict $\Ind$-normed functor $\Psi=(\Psi_A,\Psi_B):\fG(A\sqcup B)\to \fG(A)\times \fG(B)$. Let us show that $\Psi$ is an inverse of $\Phi$ in $\NSStMI$. We first construct by induction an $\Ind$-normed natural transformation $\eta:\Phi \circ \Psi \simeq \Id{\fG(A\sqcup B)}$.
		Consider any object $X=(\theta,(c_1,\ldots,c_{\theta}))$ of $\fG({A\sqcup B})$. If $\theta=\Id{}$, then $(c_1,\ldots,c_{\theta})=(c)$ and $c$ belongs either to $A$ or $B$ and we can define $\eta_{X}$ respectively by the isomorphisms $ \rho_c:c \otimes e \simeq c$ or $ \lambda_c:e \otimes c \simeq c$. Assume now that $\theta=r\bigotimes_{T}\left(\theta_1,\ldots, \theta_{|T|}\right)$, denote by $X_i$ the object $\left(\theta_i,\left(c_{|\theta_1|+\cdots +|\theta_{i-1}|}, \ldots ,c_{|\theta_1|+\cdots +|\theta_i|}\right)\right)$ and assume that $\eta_{X_i}$ is already constructed. By definition 
		$\Phi \circ \Psi (X)$ is 
		\begin{align*}
		&\Psi_A(X)\otimes \Psi_B(X) \\
		=&\textstyle r\bigotimes_T \left(\Psi_A (X_1),\ldots,\Psi_A\left(X_{|\theta|}\right)\right) \otimes r\bigotimes_T \left(\Psi_B (X_1),\ldots,\Psi_B\left(X_{|\theta|}\right)\right) \\  
		\overset{\alpha}{\simeq} &\textstyle r\bigotimes_T \left(\Psi_A (X_1)\otimes \Psi_B (X_1),\ldots, \Psi_A\left(X_{|\theta|}\right)\otimes \Psi_B \left(X_{|\theta|}\right) \right)
		\end{align*}
		and $\eta_{X}$ can be defined by $\eta_X=\bigsqcup_{i=1}^{|T|}\eta_{X_i}$ through the $G$-equivariant isomorphism $\alpha$. It follows from a direct inductive argument that $\eta$ is an $\Ind$-normed natural transformation.
We can also construct an $\Ind$-normed natural transformation $\epsilon:\Psi \circ \Phi \simeq \Id{\fG(A)\times \fG(B)}$. Given an object $X=\left(\left(\theta,\left(a_1,\ldots,a_{|\theta|}\right)\right),\left(\tau,\left(b_1,\ldots,b_{|\tau|}\right)\right)\right)$ of $\fG(A)\times \fG(B)$, $\epsilon_{X}$ is the $G$-equivariant isomorphism $\epsilon_X$ between $\left(\left(\left(\theta \otimes \tau\right)_A,(a_1,\ldots, a_{|\theta|})\right), \left(\left(\theta\otimes \tau\right)_B,\left(b_1,\ldots,b_{|\tau|}\right)\right)\right)$ and $\simeq \left(\left(\theta,\left(a_1,\ldots,a_{|\theta|}\right)\right),\left(\tau,(b_1,\ldots,b_{|\tau|}\right)\right)$ given by $\left(\Id{\num{|\theta|}},\Id{\num{|\tau|}}\right)$ where $(\theta \otimes \tau)_A$ is the formal external norm obtained by replacing each occurrence of the identity that corresponds to an element of $B$ by $\bigotimes_{\emptyset}$.
	\end{proof}
\end{proposition}

Our goal is now to use \ref{univquasi} to show that $\fG:\gset \to \NSStMI$ satisfies the $\Ind$-Beck Chevalley condition in the sense of \ref{Ichev}. For every subgroup $H\subseteq G$ and admissible map $u:A\to G/H$, we have to construct a right adjoint in $\NSStMI$ to the strict $\Ind$-normed functor $\fG(u)$. Note that in the category ${\gset}_{/(G/H)}$, the map $u:A\to G/H$ is isomorphic to the projection $\PP{T}:\ind{G}{H}T \to G/H$ where $T$ is the admissible $H$-set given by $u^{-1}\left(1_G H\right)$. We must then construct a right adjoint to $\fG(\PP{T})$ in $\NSStMI$.

\begin{construction}
	If $T$ is an admissible $H$-set, then we denote by $\fG\left(\PP{T}\right)^*:\fG(G/H)\to \fG\left(\ind{G}{H}T\right)$ the image under $\conss$ of the $G$-functor $G/H \to \fG \left(\ind{G}{H}T\right)$ that corresponds through the equivalence 
	$$
	{\func \left(G/H,\fG\left(\ind{G}{H}T\right)\right)\simeq \fG\left(\ind{G}{H}T\right)^H}
	$$ to the $H$-fixed object $\left(\bigotimes_T,\left(t_1,\ldots,t_{|T|}\right)\right)$ of $\fG\left(\ind{G}{H}T\right)$ where $T=\{t_1,\ldots,t_{|T|}\}$ is the chosen ordering on $T$.
\end{construction}

\begin{proposition}\label{GGadj}
	Let $T$ be an admissible $H$-set. The strict $\Ind$-normed functor $\fG(\PP{T})^*$ is a right adjoint of the strict $\Ind$-normed functor $\fG(\PP{T})$ in the $2$-category $\NSStMI$.
	
	\begin{proof}
		It follows from \ref{composi} that the strict $\Ind$-normed functor $\fG(\PP{T})^*\fG(\PP{T})$ is the image under $\conss$ of the $G$-functor $\ind{G}{H}T\to \fG\left(\ind{G}{H}T\right)$ that corresponds by adjunction to the $H$-functor $T\to \fG\left(\ind{G}{H}T\right)$ that sends $t$ to $\left(\bigotimes_T,\left(t_1,\ldots,t_{|T|}\right)\right)$. The unit $\eta:\Id{\fG\left(\ind{G}{H}T\right)}\to \fG(\PP{T})^*\fG(\PP{T})$ can be defined as the image under $\conss$ of the $G$-natural transformation that corresponds by adjunction to the $H$-natural transformation defined on $t_i$ by the morphism $t_i\to \left(\bigotimes_T,\left(t_1,\ldots,t_{|T|}\right)\right)$ given by the map $\num{1}\to \num{|T|}$ that sends $1$ to $i$. Using again \ref{composi}, the strict $\Ind$-normed functor $\fG(\PP{T})\fG(\PP{T})^*$ is the image under $\conss$ of the $G$-functor $G/H \to \fG(G/H)$ that corresponds through the equivalence $\func_G \left(G/H,\fG\left(\ind{G}{H}T\right)\right)\simeq \fG\left(\ind{G}{H}T\right)^H$ to the $H$-fixed object $\left(\bigotimes_T,\left(1_G,\ldots,1_G\right)\right)$ of $\fG(G/H)$ and the counit can be defined as the image under $\conss$ of the $G$-natural transformation that corresponds to the  $H$-fixed morphism $\left(\bigotimes_T,\left(1_G,\ldots,1_G\right)\right)\to 1_G$ given by the unique map $\num{|T|}\to \num{1}$. The triangle identities can be verified using \ref{composi}.
	\end{proof}
\end{proposition}

We can now show, using \ref{univquasi}, the following result.

\begin{proposition}
	The functor $\fG:\gset \to \NSStMI$ satisfies the $\Ind$-Beck-Chevalley condition in the sense of \ref{Ichev}
	\begin{proof}
		Let $u:A\to G/H$ be any admissible map. The map $u:A\to G/H$ is isomorphic in the category ${\gset}_{/(G/H)}$ to the projection $\PP{T}:\ind{G}{H}T\to G/H$ where $T$ is the finite $H$-set $u^{-1}(1_GH)$, which is admissible by \ref{adm}. It follows from \ref{GGadj} that $\fG(u)$ admits a right adjoint $\fG(u)^*$.
		Let $K\subseteq G$ be a subgroup of $G$ and let $w:G/K\to G/H$ be a $G$-equivariant map. Without loss of generality, we can assume that $K$ is included in $H$ and that $w$ corresponds to this inclusion. If we denote by $v:\ind{G}{K}\res{H}{K} T \to \ind{G}{H} T$ the $G$-equivariant map that corresponds by adjunction to the $K$-equivariant map $\res{H}{K}T\to \ind{G}{H} T$ that sends $t$ to $(1_G H,t)$, the square 
		$$
		\begin{tikzcd}
		\ind{G}{K}\res{H}{K}T \arrow[r,"v"] \arrow[d,"\PP{\res{H}{K}T}"'] & \ind{G}{H} T \arrow[d,"\PP{T}"] \\
		G/K \arrow[r,"w"'] &G/H
		\end{tikzcd}
		$$
		is Cartesian. A direct inspection shows that the mate depicted in the following diagram
		$$
		\begin{tikzcd}
		\fG\left(\ind{G}{K}\res{H}{K}T\right) \arrow[dr,Rightarrow,shorten >=6mm,shorten <=6mm] \arrow[r,"\fG(v)"]  & \fG\left(\ind{G}{H} T\right)  \\
		\fG(G/K) \arrow[r,"\fG(w)"']\arrow[u,"\fG\left(\PP{\res{H}{K}T}\right)^*"] &\fG(G/H) \arrow[u,"\fG\left(\PP{T}\right)^*"'] 
		\end{tikzcd}
		$$
		is the $\Ind$-normed natural transformation $\fG(v)\fG\left(\PP{\res{H}{K}T}\right)^*\to \fG\left(\PP{T}\right)^*\fG(w)$ given by the image under $\conss$ of the $G$-natural transformation that corresponds through the equivalence $\func_G \left(G/H,\fG\left(\ind{G}{H}T\right)\right)\simeq \fG\left(\ind{G}{H}T\right)^H$ to the $K$-fixed morphism $\left(\bigotimes_{\res{H}{K}T},\left(t_1,\ldots,t_{|T|}\right)\right)\to \left(\bigotimes_{T},\left(t_1,\ldots,t_{|T|}\right)\right)$ given by the identity $\Id{\num{|T|}}:\num{|T|}\to \num{|T|}$, which is clearly an isomorphism. The result now follows from \ref{papro}, \ref{2sa}, and \ref{univquasi}.
	\end{proof}
\end{proposition}

It follows from \ref{univprop} that $\fG:\gset \to \NSStMI$ induces a semi-additive functor $\cG : \SpaG \to N_{\Delta}\left(\NSStMI_{(2,1)}\right)$. By the following proposition, $\fG:\gset \to \NSSMI$ also satisfies the $\Ind$-Beck-Chevalley condition and $\cG$ factors as a functor $\SpaG \to N_{\Delta}\left(\NSSMI_{(2,1)}\right)$.

\begin{proposition}
	For every admissible morphism $u:A\to B$, the strict $\Ind$-normed functor $\fG(u):\fG(A)\to \fG(B)$ admits a right adjoint in the $2$-category $\NSSMI$.
	\begin{proof}
		Choose a decomposition $B\simeq \bigsqcup_{i=1}^n G/H_i$ and write $u$ through this decomposition as $\bigsqcup_{i=1}^n \PP{T_i}$ with $\PP{T_i}:\ind{G}{H_i}T_i \to G/H_i$ the canonical projection associated with an admissible finite $H_i$-set $T_i$. In the $2$-category $\NSStMI$, the strict $\Ind$-normed functor $\fG(u)$ admits a right adjoint $\fG(u)^*$ given by the strong $\Ind$-normed functor  $\fG(\PP{T_1})^*\otimes \cdots \otimes \fG\left(\PP{T_n}\right)^*$.
		Let us now show that the strong $\Ind$-normed functor $\fG(u)^*$ is isomorphic to a strict $\Ind$-normed functor.
		Every $\fG(\PP{T_i})^*$ is a strict $\Ind$-normed functor provided as the image of some $G$-functor under $\conss$ and it is enough to show, given two $G$-functors $\phi,\psi:A\rightrightarrows \fG(B)$, that the strong $\Ind$-normed functor $\conss(\phi)\otimes \conss (\psi)$ is isomorphic to a strict $\Ind$-normed functor. More precisely, a direct induction shows that the strong $\Ind$-normed functor $\conss(\phi)\otimes \conss (\psi)$ is isomorphic to the image under $\conss$ of the $G$-functor given by the composition $A \overset{(\phi,\psi)}{\longrightarrow}\fG (B)\times \fG(B) \overset{\otimes}{\longrightarrow}\fG(B)$. The result follows.
	\end{proof}
\end{proposition}

\subsection{The proof of the main Theorem}\label{subsecmainth}
In this section, we prove \ref{THETH}. 
We have morphisms of simplicial sets 
$$
\SpaG \overset{\cG}{\longrightarrow} N_{\Delta} \left(\NSSMI_{(2,1)}\right) \overset{N_{\Delta}(\iota)}{\longrightarrow}N_{\Delta} \left(\NSSMI \right) \overset{}{\hookrightarrow} N_{\Delta} \left(\Alg_{\EII}\left(\gsset\right)\right)
$$
where $\iota:\NSSMI_{(2,1)} \to \NSSMI$ denotes the strict $2$-functor that sends an $\Ind$-normed symmetric monoidal category to its core and where the third morphism is given by \ref{subalgNS}. We show now that the composition of these morphisms factors as a functor $\theta:\SpaG \to \Ninf$. This composition sends a finite $G$-set $A$ to the nerve of the core of the $\Ind$-normed symmetric monoidal category $\fG(A)$ and we have to show that this object is fibrant and cofibrant in $\Alg_{\EII}\left(\gsset \right)$. If we denote by $\sF(A)$ the core of $\fG(A)$, the fibrancy of $N(\sF(A))$ follows from the fact that its underlying $G$-simplicial set is the nerve of a $G$-groupoid. To show that $N(\sF(A))$ is cofibrant, we show that it is the free $\EII$-algebra over $A$. We show first that $\sF(A)$ is the free $\Ind$-normed symmetric monoidal category over $A$ seen as a discrete $G$-category. We start with the following elementary lemma.

\begin{lemma}\label{lemtec}
Let $\Gamma$ be a discrete group and let $X$ and $Y$ be two $G\times \Gamma$-sets such that the $\Gamma$-action on $X$ is free. The canonical $G$-equivariant map 
$$
f\colon \left(X^{n+1}\times Y\right)_{\Gamma} \to \underbrace{\left(X\times X \times Y\right)_\Gamma\times_{(X\times Y)_\Gamma}\cdots\times_{(X\times Y)_\Gamma} (X\times X \times Y)_\Gamma}_{n \text{ times}}
$$
is a bijection.
\begin{proof}
We prove first that $f$ is injective. Consider two tuples $A=(x_1,\ldots,x_{n+1},y)$ and $B=(x_1',\ldots,x_{n+1}',y')$ in $\left(X^{n+1}\times Y\right)$ with the same image under $f$. This means there are elements $\gamma_1,\ldots,\gamma_n$ in $\Gamma$ such that for all $i=1,\ldots,n$, we have 
$$
\gamma_i  (x_i,x_{i+1},y)=(x_i',x_{i+1}',y').
$$
In particular, it follows that $\gamma_i x_{i+1} =x_{i+1}'$ and $\gamma_{i+1}x_{i+1}=x_{i+1}'$ and we deduce that $\gamma_i=\gamma_{i+1}$ using that the $\Gamma$-action on $X$ is free. This implies that $A$ and $B$ represent the same element of $\left(X^{n+1}\times Y\right)_{\Gamma}$. We prove now that $f$ is surjective. Consider an element $A$ of 
$$
\underbrace{\left(X\times X \times Y\right)_\Gamma\times_{(X\times Y)_\Gamma}\cdots\times_{(X\times Y)_\Gamma} (X\times X \times Y)_\Gamma}_{n \text{ times}}
$$
represented by tuples $(x_i,z_i,y_i)$ for all $i=1,\ldots,n$. By assumptions, there exists $\gamma_i$ in $\Gamma$ such that $z_i=\gamma_i x_{i+1}$ and $\gamma_i y_{i+1}=y_i$ for all $i=1,\ldots,n-1$. For all $i=1,\ldots,n+1$, consider the element
$$
w_i=\begin{cases}
x_1 & \text{if} ~~ i=1; \\
z_1 & \text{if} ~~ i=2; \\
(\gamma_1\gamma_2\cdots\gamma_{i-2})\cdot  z_{i-1} & \text{if} ~~ i=3,\ldots,n+1.
\end{cases}
$$
Note that we have $(w_i, w_{i+1}, y_1 )=(\gamma_1\gamma_2 \cdots \gamma_{i-1} )\cdot (x_i,z_i,y_i)$ and, in particular, the image of the orbit of $(w_1,\ldots, w_{n+1},y_1)$ under $f$ is $A$.

\end{proof}
\end{lemma}
One important direct consequence of \ref{lemtec} is the following result.
\begin{proposition}\label{propquogcat}
Let $\Gamma$, $X$ and $Y$ be as in \ref{lemtec}. The $G$-category $\left(\widetilde{X}\times Y\right)_{\Gamma}$ is the $G$-category whose $G$-set of objects is $(X\times Y)_{\Gamma}$, whose $G$-set of morphisms is $(X\times X \times Y)_{\Gamma}$ and whose composition is given by the $G$-equivariant map
$$
(X\times X \times Y)_{\Gamma}\times_{(X\times Y)_{\Gamma}}(X\times X \times Y)_{\Gamma} \simeq (X\times X\times X \times Y)_{\Gamma}\overset{((p_1,p_3)\times  \Id{Y})_{\Gamma} }{\longrightarrow}
(X\times X \times Y)_{\Gamma}.
$$
\end{proposition}

Using this proposition, one can prove that $\sF (A)$ is the free $\Ind$-normed symmetric monoidal category on $A$.

\begin{proposition}\label{free}
	The $\Ind$-normed symmetric monoidal category $\sF (A)$ is free over $A$ for strict $\Ind$-normed functors. Moreover, for every $\Ind$-normed symmetric monoidal category $\cC$, the $G$-functor $A\to \sF(A)$ induces an equivalence of categories 
	$$
	\NSSMI (\sF (A),\cC)\simeq \func_G (A,\cC).
	$$
	\begin{proof}
		By \ref{rubco}, it is enough to show that $\sF(A)$ is the $\EI$-algebra given by the image of $A$ under the monad $\EI:\catG \to \catG$. Given a $G$-category $\cC$, if we denote respectively by $\Ob(\cC)$ and $\Mor(\cC)$ the associated $G$-sets of objects and morphisms, then, using \ref{propquogcat}, we have
		$$
		\Ob \left(\bigsqcup_{n\ge 0}(\EI)_n\times_{\Sigma_n}A^n\right) \simeq   \bigsqcup_{n\ge 0} \SN_n \times_{\Sigma_n} A^n\simeq \bigsqcup_{n\ge 0} F(V_{\cI})_n\times A^n \simeq \Ob (\sF(A))
		$$
        and
        \begin{align*}
		\Mor \left(\bigsqcup_{n\ge 0}(\EI)_n\times_{\Sigma_n}A^n\right) &\simeq \bigsqcup_{n\ge 0} {\SN_n}^2 \times_{\Sigma_n} A^n\simeq \bigsqcup_{n\ge 0} {F(V_{\Ind})_n}^2\times A^n \times \Sigma_n
        \\
         & \simeq \Mor (\sF(A)),
		\end{align*}
		with $F(V_{\Ind})$ defined in \ref{rubmod}.
        It is clear that these isomorphisms induce an isomorphism $\sF (A)\simeq \EI (A)$. 
		The equivalence of categories 
		$$
		\NSSMI (\FF (A),\cC)\simeq \func_G (A,\cC)
		$$
		is a consequence of the freeness of $\sF(A)$ and of the description of $\NSSMI (\sF (A),\cC)$ in terms of the arrow category of $\cC$ given by \ref{aromor}.
	\end{proof}
\end{proposition}

\begin{remark}
	Given a $G$-functor $\phi:A\to \fG(B)$, the core of the strict $\Ind$-normed functor $\conss (\phi)$ is the strict $\Ind$-normed functor $\sF(A)\to \sF(B)$ induced by $\phi:A\to \sF(B)$ through the equivalence $\NSSMI (\sF(A),\sF(B))\simeq \func_G(A,\sF(B))$ of \ref{free}.
\end{remark}

The following lemma is the key point to transfer the freeness of $\sF(A)$ as an $\Ind$-normed symmetric monoidal category to the freeness of $N(\sF(A))$ as an $\EII$-algebra.

\begin{lemma}
	Let $\Gamma$, $X$ and $Y$ be as in \ref{lemtec}. The natural $G$-equivariant map
	$$
	N\left(\widetilde{X}\times Y\right)/\Gamma \to N\left(\left(\widetilde{X}\times Y\right)/\Gamma\right)
	$$
	is an isomorphism of $G$-simplicial sets where $\widetilde{X}$ is the translation $G$-groupoid of $X$, as defined in \ref{transg}.
	\begin{proof}
	This is a direct consequence of \ref{lemtec}.
	\end{proof}
\end{lemma}

\begin{corollary}
	Let $A$ be a finite $G$-set. The nerve of $\sF(A)$ is the free $N_{\infty}$-algebra over $A$.
	\begin{proof}
		We have to show that the nerve of $\sF(A)$ is the image of $A$ by the monad $\EII:\gsset \to \gsset$. We proved in \ref{free} that $\sF(A)$ is the free $\Ind$-normed symmetric monoidal category over $A$ and by definition $(\EII)_n=N((\EI)_n)$. It is then enough to show that the canonical morphism
		$$
		\bigsqcup_{n\ge 0}A^n\times_{\Sigma_n} N\left((\EI)_n\right)\to N\left(\bigsqcup_{n \ge 0} A^n\times_{\Sigma_n} (\EI)_n\right)
		$$
		is an isomorphism. This is a consequence of the fact that the nerve preserves finite coproducts and of the previous lemma using that $(\EI)_n=\widetilde{\SN_n}$.
	\end{proof}
\end{corollary}
If $A$ is a finite $G$-set, then the $\EII$-algebra $N(\sF(A))$ is fibrant and cofibrant and we get a functor $\theta:\SpaG \to \Ninf$ that sends $A$ to $N(\sF(A))$. We will show now that this functor is fully faithful. For every subgroup $H\subseteq G$ and finite $G$-set $A$, there are equivalences 
$$
\NSSMI (\sF(G/H),\sF(A))\simeq \func_G (G/H,\sF(A)) \simeq \sF(A)^H
$$
and the latter groupoid can be described by the following proposition.
\begin{proposition}\label{fixed}
	For every finite $G$-set $A$, there is an equivalence of categories
	$$
	\Ind (H)_{/\res{G}{H}A} \simeq \fG(A)^H.
	$$
	\begin{proof}
		Let $u:T \to \res{G}{H}A$ be an $H$-equivariant map with $T$ an admissible $H$-set. Consider the $H$-fixed object $\left(\bigotimes_T,\left(u(t_1),\ldots,u\left(t_{|T|}\right)\right)\right)$ of $\fG(A)$ where $T= \{t_1,\ldots,t_{|T|}\}$ is given by the chosen ordering on $T$. If $w:u\to v$ is a morphism in $\Ind (H)_{/\res{G}{H}A}$ between two $H$-equivariant maps $u:T\to \res{G}{H}A$ and $v:S\to \res{G}{H}A$, we get an $H$-fixed morphism $\left(\bigotimes_T,\left(u(t_1),\ldots,u\left(t_{|T|}\right)\right)\right) \to \left(\bigotimes_S,\left(v(s_1),\ldots, v\left(s_{|S|}\right)\right)\right)$ given by the map $w:\num{|T|}\to \num{|S|}$ that corresponds to $w:T\to S$ through the chosen orderings $T= \{t_1,\ldots, t_{|T|}\}$ and $S= \{s_1,\ldots, s_{|S|}\}$. We thus get a functor $\Ind (H)_{/\res{G}{H}A}\to \fG(A)^H$ which is fully faithful. Given now an $H$-fixed object $\left(\theta, \left(a_1,\ldots,a_{|\theta|}\right)\right)$ of $\fG(A)^H$, the identity $\Id{\num{|\theta|}}:\num{|\theta|}\to \num{|\theta|}$ induces an $H$-fixed isomorphism $\left(\bigotimes_{T_\theta},\left(u(t_1),\ldots, u\left(t_{|\theta|}\right) \right)\right)\to \left(\theta, \left(a_1,\ldots,a_{|\theta|}\right)\right)$ where $u:T_{\theta}\to A$ is the $H$-equivariant map that sends $t_i$ to $a_i$. We conclude that the functor $\Ind (H)_{/\res{G}{H}A} \to \fG(A)^H$ is essentially surjective and the result follows.
	\end{proof}
\end{proposition}

The $\Ind$-normed symmetric monoidal category $\sF(A)$ is by definition the core of $\fG(A)$ and we deduce the following corollary.

\begin{corollary}\label{hfixob}
	For every finite $G$-set $A$, there is an equivalence of groupoids 
	$$
	\iota \left(\Ind (H)_{/\res{G}{H}A} \right)\simeq \sF(A)^H.
	$$
\end{corollary}

We can now use these descriptions to show the following result.

\begin{proposition}\label{equi21}
	The functor $\theta:\SpaG \to \Ninf$ is fully faithful.
	\begin{proof}
		By definition, the functor $\theta:\SpaG \to \Ninf$ is the composition of morphisms of simplicial sets 
		$$
		\SpaG \overset{\cG}{\longrightarrow} N_{\Delta} \left(\NSSMI_{(2,1)}\right) \overset{N_{\Delta}(\iota)}{\longrightarrow}N_{\Delta} (\NSSMI) \overset{}{\hookrightarrow} N_{\Delta} \left(\Alg_{\EII}\left(\gsset\right)\right)
		$$
		which factors through the inclusion $\Ninf \hookrightarrow N_{\Delta} \left(\Alg_{\EII}\left(\gsset\right)\right)$. By \ref{thesubOalg}, the simplicial functor $ \NSSMI \overset{}{\hookrightarrow} \Alg_{\EII}\left(\gsset\right)$ is fully faithful and it is enough to show that the functor $N_{\Delta}(\iota) \circ \cG:\SpaG \to N_{\Delta} (\NSSMI)$ is fully faithful.
		The $\infty$-category $\SpaG$ is a $(2,1)$-category and, if $A$ and $B$ are finite $G$-sets, the groupoid $\Hom_{\SpaG}(A,B)$ can be described as follow: the objects in $\Hom_{\SpaG}(A,B)$ are spans $A \overset{f}{\longleftarrow} C \overset{g}{\longrightarrow} B$ and morphisms are isomorphisms of spans
		$$
		\begin{tikzcd}[row sep=small]
		&C \arrow[rd,"g"] \arrow[dd,"\simeq"] \arrow[ld,"f"']& \\
		A& & B. \\
		&D \arrow[ur,"i"'] \arrow[ul,"h"]&
		\end{tikzcd}
		$$
		Seen as a functor $\SpaG \to N_{\Delta}( \NSStMI )$, the functor $N_{\Delta}(\iota) \circ \cG$ is semi-additive by \ref{univquasi} and it is enough to show, given a subgroup $H\subseteq G$ and a finite $G$-set $A$, that $(N_{\Delta}(\iota) \circ \cG)_{G/H,A}:\Hom_{\SpaG}(G/H,A)\to \Hom_{N_{\Delta}(\NSSMI)}(\sF(G/H),\sF(A))$ is an equivalence. Every span in $\Hom_{\SpaG}(G/H,A)$ is isomorphic to a span of the form $G/H \overset{\PP{T}}{\longleftarrow} \ind{G}{H} T \overset{}{\longrightarrow} A$ with $T$ an admissible $H$-set and there is an equivalence of groupoids $\Hom_{\SpaG}(G/H,A) \simeq \iota \left(\Ind (H)_{/\res{G}{H}A} \right)$ that sends an $H$-equivariant map $u:T\to A$ to the span $G/H \overset{}{\longleftarrow} \ind{G}{H} T \overset{v}{\longrightarrow} A$ with $v:\ind{G}{H}T \to A$ the $G$-equivariant map that corresponds to $u$ by adjunction. By definition, the functor $(N_{\Delta}(\iota) \circ \cG)_{G/H,A}:\Hom_{\SpaG}(G/H,A)\to \Hom_{N_{\Delta}(\NSSMI)}(\sF(G/H),\sF(A))$ sends the span $G/H \overset{}{\longleftarrow} \ind{G}{H} T \overset{v}{\longrightarrow} A$ to the composition $\fG(v)\fG(\PP{T})^*$ which is the image under $\conss$ of the $G$-functor $G/H \to \sF(A)$ that corresponds to the $H$-fixed object $\left(\bigotimes_T,\left(u(t_1),\ldots,u\left(t_{|T|}\right)\right)\right)$. Moreover, the functor $(N_{\Delta}(\iota) \circ \cG)_{G/H,A}:\Hom_{\SpaG}(G/H,A)\to \Hom_{N_{\Delta}(\NSSMI)}(\sF(G/H),\sF(A))$ sends an isomorphism of spans
		$$
		\begin{tikzcd}[row sep=small]
		&\ind{G}{H}T \arrow[rd,"u"] \arrow[dd,"\ind{G}{H} w"] \arrow[ld,"\PP{T}"']& \\
		G/H& & A \\
		&\ind{G}{H}S \arrow[ur,"v"'] \arrow[ul,"\PP{S}"]&
		\end{tikzcd}
		$$
		to the $\Ind$-normed natural transformation
		$$
		\fG(u) \fG(\PP{T})^* \simeq \fG(v) \fG\left(\ind{G}{H}w\right) \fG\left(\ind{G}{H}w^{-1}\right) \fG(\PP{S})^* =\fG(v)\fG(\PP{S})^*
		$$
		which corresponds through the equivalence $\NSSMI (\sF(G/H),\sF(A))\simeq \sF(A)^H$ to the isomorphism 
		$$\textstyle
		\bigotimes_T\left(u(t_1),\ldots,u\left(t_{|T|}\right)\right)\to \bigotimes_S \left(v(s_1),\ldots,v\left(s_{|S|}\right)\right)
		$$
		given by $w:T\to S$ seen as a map $w:\num{|T|}\to \num{|S|}$. 
		Finally, it follows that the functor $(N_{\Delta}(\iota) \circ \cG)_{G/H,A}:\Hom_{\SpaG}(G/H,A)\to \Hom_{N_{\Delta}(\NSSMI)}(\sF(G/H),\sF(A))$ corresponds to the equivalences 
		$$
		\Hom_{\SpaG}(G/H,A)\simeq \iota \left(\Ind (H)_{/\res{G}{H}A} \right) \simeq \sF(A)^H \simeq 
		\NSSMI (\sF(G/H),\sF(A)).
		$$
	\end{proof}
\end{proposition}

We can finally prove our main theorem.
\begin{theorem}\label{THETHH}
	The left Kan extension $\Theta:\Mack \to \Ninf$ of $\theta:\SpaG \to \Ninf$ along the Yoneda embedding $\SpaG \to \Mack$ is an equivalence of $\infty$-categories.
	\begin{proof}
        
	    Combining \cite[Corollary 3.22]{greg3} with the Elmendorf's Theorem, we obtain that the free-forgetful adjunction 
		$$
		F:\gs \rightleftarrows \Ninf :U
		$$
		satisfies the conditions of \cite[Corollary 4.7.3.18]{HA}. It follows that $\Ninf$ is projectively generated and that the $\EII$-algebras $N(\sF(A))$ for all finite $G$-sets $A$ form a set of compact projective generators. Finally, it follows from \ref{equi21} and \cite[Proposition 5.5.8.22]{MR2522659} that the left Kan extension $\Theta:\Mack \to \Ninf$ of $\theta:\SpaG \to \Ninf$ is an equivalence of $\infty$-categories.
	\end{proof}
\end{theorem}

\begin{remark}
	By definition, the inverse of the functor $\Theta:\Mack \to \Ninf$ sends a (fibrant) $N_{\infty}$-algebra $X$ to the $\Ind$-Mackey functor $\ma{M}_X:\opo{\SpaG} \to \sS$ that sends $G/H$ to the space 
	$$
	\Hom_{\Ninf}(N(\sF(G/H)),X)\simeq X^H.
	$$
	If $K\subseteq H$ are subgroups of $G$, then the functor $\res{H}{K}:X^H \to X^K$ obtained as the image of the span $G/K \overset{\Id{G/K}}{\longleftarrow} G/K \overset{u}{\longrightarrow} G/H$ by $\ma{M}_X$ corresponds to the canonical inclusion $X^H \hookrightarrow X^K$. If $H/K$ is admissible, then the functor $\TTT_{K}^H:X^K \to X^H$ obtained as the image of the span $G/H \overset{u}{\longleftarrow} G/K \overset{\Id{G/K}}{\longrightarrow} G/K$ by $\ma{M}_X$ sends a $n$-simplex $x$ of $X^K$ to the $n$-simplex $\bigotimes_{H/K}^X\left(h_1x,\ldots ,h_{[H:K]}x\right)$ of $X^H$ where $h_1, \ldots, h_{[H:K]}$ are representatives of $H/K$ and with $\bigotimes_{H/K}^X:X^{\times [H:K]} \to X$ the image of $\bigotimes_{H/K}$ by the morphism of operads $|-|_X:\EII \to \End_X$ that corresponds to the $\EII$-algebra structure of $X$.
\end{remark}

\section{Appendix: Operads in \texorpdfstring{$G$}{G}-categories}\label{ApA}

In this appendix, we recall the definitions of the $2$-categories of $\ma{O}$-algebras over an operad $\ma{O}$ in $\Cat^G$. The goal is to prove \ref{thesubOalg}, which is used to prove \ref{equi21}. We follow \cite[Section $4$]{Rub2} for the definitions and refer the reader to it for more details. 

\begin{notation}
	If $\ma{O}$ is an operad in $\catG$ and $\cC$ an $\ma{O}$-algebra in $\catG$, then we will denote by $|-|_{\cC}:\ma{O}\to \End_{\cC}$ the morphism of operads in $\catG$ that corresponds to the structure of $\ma{O}$-algebra of $\cC$.
\end{notation}

\begin{definition}[\protect{\cite[Definition 4.2]{Rub2}}]
	Let $\ma{O}$ be an operad in $\catG$ and $\cC$ and $\cD$ two $\ma{O}$-algebras in $\catG$. A \emph{lax $\ma{O}$-algebra morphism $\phi:\cC \to \cD$} consists of:
	\begin{enumerate}
		\item a $G$-functor $\phi :\cC \to \cD$;
		\item for every $n\ge 0$ and $x \in \ma{O}_n$, a natural transformation $\phi _x:|x|_{\cD} \circ \phi^{\times n}\Rightarrow \phi \circ |x|_{\cC}$ of functors $\cC^{\times n} \to \cD$;
	\end{enumerate}
	such that the following conditions are satisfied:
	\begin{enumerate}
		\item for every $n\ge 0$, the maps $\phi_x$ vary naturally in $x\in \ma{O}_n$;
		\item for every $n\ge 0$, $x\in \ma{O}_n$ and $(g,\sigma)\in G\times \Sigma_n$, the equation $(g,\sigma)\cdot \phi_x=\phi_{(g,\sigma)\cdot x}$ holds, ie, $\phi_x$ is $(G\times \Sigma_n)$-equivariant;
		\item $\phi_{\Id{\ma{O}}}=\Id{\phi}$;
		\item for every $y\in \ma{O}_n$ and $x_i\in \ma{O}_{k_i}$, the  transformation $\phi_{\gamma (y,x_1,\ldots,x_n)}$ is equal to the composite 
		$$
		\begin{tikzcd}
		{|y|}_{\cD}\circ \left( \left({|x_1|}_{\cD}\circ \phi^{\times k_1}\right)\times \cdots \times \left({|x_n|}_{\cD}\circ \phi^{\times k_n}\right)\right) \arrow[d,Rightarrow,"\Id{|y|_{\cD}}\circ \left(\phi_{x_1}\times \cdots\times \phi_{x_n}\right)"] \\
		\left({|y|}_{\cD}\circ \phi^{\times n}\right)\circ \left({|x_1|}_{\cC}\times \cdots \times {|x_n|}_{\cC}\right) \arrow[d,Rightarrow,"\phi_y \circ \left(\Id{|x_1|_{\cC}}\times \cdots \times \Id{|x_n|_{\cC}} \right)"] \\
		\phi \circ {|y|}_{\cC}\circ \left({|x_1|}_{\cC}\times \cdots \times {|x_n|}_\cC\right).
		\end{tikzcd}
		$$
	\end{enumerate}
	A \emph{strong $\ma{O}$-algebra morphism (resp. strict $\ma{O}$-algebra morphism)} is a lax $\ma{O}$-algebra morphism $\phi:\cC \to \cD$ such that $\phi_x$ is an isomorphism (resp. the identity) for every $n\ge 0$ and $x\in \ma{O}_n$.
\end{definition}

If $\phi:\cC \to \cD$ and $\psi:\cD \to \ma{E}$ are two lax $\ma{O}$-algebra morphisms, the composition $\psi \circ \phi$ has a natural structure of lax $\ma{O}$-algebra morphism where for every $n\ge 0$ and $x\in \ma{O}_n$, the natural transformation $(\psi \circ \phi)_x$ is given by the composition 
$$
|x|_{\ma{E}}\circ (\psi \circ \phi)^{\times n}=|x|_{\ma{E}}\circ \psi^{\times n}\circ \phi^{\times n} \overset{\psi_x \circ \left(\Id{\phi^{\times n}}\right) }{\longrightarrow} \psi \circ |x|_{\cD} \circ \phi^{\times n} \overset{\Id{\psi}\circ \phi_x}{\longrightarrow} \psi \circ \phi \circ |x|_{\cC}.
$$
If $\phi$ and $\psi$ are strong $\ma{O}$-algebra morphisms (resp. strict), their composition $\psi  \circ \phi$ is also a strong $\ma{O}$-algebra morphism (resp. strict).
\begin{definition}[\protect{\cite[Definition 4.3]{Rub2}}]
	Let $\ma{O}$ be an operad in $\catG$ and $\phi,\psi:\cC \rightrightarrows \cD$ a pair of lax $\ma{O}$-algebra morphisms.
	An \emph{$\ma{O}$-natural transformation} $\omega:\phi \Rightarrow \psi$ is a $G$-natural transformation $\omega:\phi \Rightarrow \psi$ such that for every $n\ge 0$ and $x\in \ma{O}_n$, the following diagram of natural transformations commutes
	$$
	\begin{tikzcd}
	{|x|}_{\cD}\circ \phi^{\times n} \arrow[r,Rightarrow,"\Id{|x|_{\cD}}\circ \omega^{\times n}"] \arrow[d,Rightarrow,"\phi_x"'] & {|x|}_{\cD}\circ \psi^{\times n}  \arrow[d,Rightarrow,"\psi_x"] \\ 
	\phi \circ |x|_{\cC} \arrow[r,Rightarrow,"\omega \circ \Id{{|x|}_{\cC}} "'] & \psi\circ |x|_{\cC}.
	\end{tikzcd}
	$$
\end{definition}

\begin{notation}\label{2catalg}
	We denote by $\OAlgL$ (resp. $\OAlgStg$ and $\OAlgSt$) the $2$-category of $\ma{O}$-algebras in $\catG$, lax $\ma{O}$-algebra morphisms (resp. strong and strict) and $\ma{O}$-natural transformations.
\end{notation}

Let $\ma{O}$ be an operad in $\catG$, $\cC$ an $\ma{O}$-algebra in $\catG$ and $\cD$ a $G$-category. The functor $\func (\cD, (-)):\catG \to \catG$ preserves products and $\func (\cD,\cC)$ is naturally a $\func (\cD,\ma{O})$-algebra. It follows that the $G$-category of nonequivariant functors $\func (\cD,\cC)$ is naturally endowed with the structure of an $\ma{O}$-algebra where the morphism of $G$-operads $|-|_{\func (\cD,\cC)}:\ma{O}\to \End_{\func (\cD,\cC)}$ is given by the composition 
$$
\ma{O} \to \func (\cD,\ma{O}) \overset{|-|_{\func(\cC,\cD)}}{\longrightarrow} \End_{\func (\cD,\cC)}.
$$
Here the first morphism $\ma{O} \simeq \func (*,\ma{O}) \to \func (\cD,\ma{O})$ is induced by the unique functor $\cD \to *$.
In particular, the arrow category $\Arr(\cC)=\func (\Delta^1,\cC)$ of $\cC$ has a natural structure of an $\ma{O}$-algebra.

\begin{remark}
	If $\ma{O}$ is an operad in $\catG$, $\cC$ an $\ma{O}$-algebra in $\catG$ and $\cD$ a $G$-category, then, for every $n\ge 0$ and $x\in \ma{O}_n$, the operation $|x|_{\func (\cD,\cC)}$ is given by the composition $\func (\cD,\cC)^{\times n} \overset{\times }{\longrightarrow} \func\left(\cD^{\times n},\cC^{\times n}\right) \overset{(-)\circ |x|_{\cC}}{ \longrightarrow} \func(\cD,\cC)$.
\end{remark}

If $\cC$ and $\cD$ are two $G$-categories, the data of a $G$-natural transformation $\omega:\phi \Rightarrow \psi$ between two $G$-functors $\phi, \psi :\cC\rightrightarrows \cD$ is equivalent to the data of a $G$-functor $\Omega:\cC \to \Arr (\cD)$. This characterisation of $G$-natural transformations lifts to the case of $\ma{O}$-natural transformations.  
\begin{proposition}
	Let $\cC$ and $\cD$ be two $\ma{O}$-algebras in $\CatG$. The data of an $\ma{O}$-natural transformation $\omega:\phi \Rightarrow \psi$ between two lax $\ma{O}$-algebra morphisms (resp. strong or strict) $\phi,\psi:\cC \rightrightarrows \cD$ is equivalent to the data of a lax $\ma{O}$-algebra morphism (resp. strong or strict) $\Omega:\cC \to \Arr (\cD)$. 
	\begin{proof}
		Let $\omega:\phi \Rightarrow \psi$ be an $\ma{O}$-natural transformation between two lax $\ma{O}$-algebra morphisms (resp. strong or strict) $\phi,\psi:\cC \rightrightarrows \cD$. The $G$-functor $\Omega:\cC \to \Arr (\cD)$ associated with the $G$-natural transformation $\omega$ admits a structure of lax $\ma{O}$-algebra morphism (resp. strong or strict)  where for $n\ge 0$ and $x\in \ma{O}_n$, the natural transformation $\Omega_x:|x|_{\Arr (\cD)} \circ \Omega^{\times n}\Rightarrow \Omega \circ |x|_{\cC}$ is given by the diagram 
		$$
		\begin{tikzcd}
		{|x|}_{\cD}\circ \phi^{\times n} \arrow[r,Rightarrow,"\Id{|x|_{\cD}}\circ \omega^{\times n}"] \arrow[d,Rightarrow,"\phi_x"'] & {|x|}_{\cD}\circ \psi^{\times n}  \arrow[d,Rightarrow,"\psi_x"] \\ 
		\phi \circ |x|_{\cC} \arrow[r,Rightarrow,"\omega \circ \Id{{|x|}_{\cC}} "'] & \psi\circ |x|_{\cC}
		\end{tikzcd}
		$$
		seen as a morphism in the category $\func (\cC^{\times n},\Arr (\cD))\simeq \Arr (\func (\cC^{\times n},\cD))$, which is commutative by assumption using that $\omega$ is an $\ma{O}$-natural transformation. It follows directly from the fact that $\phi$ and $\psi$ are lax $\ma{O}$-algebra morphisms that $\Omega$ is a lax $\ma{O}$-algebra morphism (resp. strong or strict). Assume now that $\Omega:\cC \to \Arr (\cD)$ is a lax $\ma{O}$-algebra morphism (resp. strong or strict) and consider the $G$-natural transformation $\omega:\phi \Rightarrow \psi$ associated with $\Omega$ seen as a $G$-functor. The $G$-functors $\phi$ and $\psi$ are by definition respectively defined by $d^0\circ \Omega$ and $d^1\circ \Omega$ where $d^0:\Arr(\cD) \to \cD$ and $d^1:\Arr (\cD) \to \cD$ are the strict $\ma{O}$-algebra morphisms induced by the inclusions $\partial^0:\Delta^0 \to \Delta^1$ and $\partial^1:\Delta^0 \to \Delta^1$. It follows that $\phi$ and $\psi$ are naturally lax $\ma{O}$-algebra morphisms (resp. strong and strict). Finally, the natural transformation $\Omega_x:|x|_{\Arr (\cD)} \circ \Omega^{\times n}\Rightarrow \Omega \circ |x|_{\cC}$ seen as a morphism in $\Arr (\func (\cC^{\times n},\cD)) \simeq \func (\cC^{\times n},\Arr (\cD))$ directly ensures that the diagram 
		$$
		\begin{tikzcd}
		{|x|}_{\cD}\circ \phi^{\times n} \arrow[r,Rightarrow,"\Id{|x|_{\cD}}\circ \omega^{\times n}"] \arrow[d,Rightarrow,"\phi_x"'] & {|x|}_{\cD}\circ \psi^{\times n}  \arrow[d,Rightarrow,"\psi_x"] \\ 
		\phi \circ |x|_{\cC} \arrow[r,Rightarrow,"\omega \circ \Id{{|x|}_{\cC}} "'] & \psi\circ |x|_{\cC}
		\end{tikzcd}
		$$
		is commutative and it follows that $\omega$ is an $\ma{O}$-natural transformation. The result follows.
	\end{proof}
\end{proposition}
If $\ma{O}$ is an operad in $\catG$ and $\cC$ and $\cD$ two $\ma{O}$-algebras in $\CatG$, then it follows directly from the previous proposition that the categories $\OAlgL (\cC,\cD)$, $\OAlgStg (\cC,\cD)$ and $\OAlgSt (\cC, \cD)$ admit the following description.

\begin{proposition}\label{aromor}
	Let $\ma{O}$ be an operad in $\catG$ and $\cC$ and $\cD$ two $\ma{O}$-algebras in $\CatG$. If $\OAlg (\cC,\cD)$ denotes either $\OAlgL(\cC,\cD)$, $\OAlgStg(\cC,\cD)$ or $\OAlgSt (\cC,\cD)$ and if $\Oalg(\cC,\cD)=\Ob (\OAlg (\cC,\cD))$
	denotes the corresponding class of $G$-functors, then the category given by the following diagram 
	$$
	\begin{tikzcd}[font=\small,column sep=small]
	\Oalg (\cC,\Arr (\cD))\times_{\Oalg(\cC,\cD)} \Oalg (\cC,\Arr (\cD)) \arrow[r]& \Oalg(\cC,\Arr(\cD)) \arrow[r,shift right=0.2cm] \arrow[r,shift left=0.2cm] & \arrow[l] \Oalg (\cC,\cD)
	\end{tikzcd}
	$$
	induced by the canonical strict $\ma{O}$-algebra morphisms $\Arr (\cD)\times_{\cD}\Arr (\cD) \to \Arr (\cD)$, $\Arr(\cD) \rightrightarrows \cD$ and $\cD\to \Arr (\cD)$ is equivalent to $\OAlg (\cC,\cD)$.
\end{proposition}
Finally, we prove the important result of this appendix.

\begin{proposition}\label{thesubOalg}
	There is a fully faithful simplicial functor $N\colon \OAlgSt \to \Alg_{N(\ma{O})}\left(\gsset\right)$ that sends an $\ma{O}$-algebra $\cC$ in $\CatG$ to its nerve $N(\cC)$ where the strict $2$-category $\OAlgSt$ is seen as a simplicial category.
	\begin{proof}
		Let $\cC$ and $\cD$ be two $\ma{O}$-algebras in $\CatG$. It follows from \ref{aromor} that we have bijections
		\begin{align*}
		& N(\OAlgSt(\cC,\cD))_n \\
		= & \Oalg (\cC,\Arr (\cD))\times_{\Oalg(\cC,\cD)} \cdots \times_{\Oalg(\cC,\cD)} \Oalg (\cC,\Arr (\cD))  \\
		\simeq & \Oalg (\cC,\Arr (\cD)\times_{\cD} \cdots \times_{\cD} \Arr(\cD)) \\
		\simeq & \Alg_{N(\ma{O})}\left(\gsset\right) \left(N(\cC),N(\cD)^{\Delta^n}\right) \\
		=& \Hom_{\Alg_{N(\ma{O})}\left(\gsset\right)}(N(\cC),N(\cD))_n
		\end{align*}
		which give an isomorphism of simplicial sets 
        $$
        N(\OAlgSt(\cC,\cD))\simeq \Hom_{\Alg_{N(\ma{O})}\left(\gsset\right)}(N(\cC),N(\cD)).
        $$
		It is clear that these isomorphisms induce a fully faithful simplicial functor $\OAlgSt \to \Alg_{N(\ma{O})}\left(\gsset\right)$.
	\end{proof}
\end{proposition}

\bibliographystyle{amsalpha}
\bibliography{references}

\end{document}